\documentclass{amsart}
\usepackage{amssymb,amsfonts}
\usepackage{enumerate}
\usepackage{enumitem}
\usepackage{mathrsfs}
\usepackage{graphicx}
\usepackage{mathtools}
\usepackage{leftidx}
\usepackage{tikz-cd}
\usepackage{amsmath,calligra,mathrsfs}
\usepackage{bm}
\usepackage{url}
\newtheorem{thm}{Theorem}[section]
\newtheorem{cor}[thm]{Corollary}
\newtheorem{prop}[thm]{Proposition}
\newtheorem{lem}[thm]{Lemma}

\theoremstyle{definition}
\newtheorem{defn}[thm]{Definition}

\newtheorem{exmp}[thm]{Example}

\newtheorem{conv}[thm]{Convention}

\theoremstyle{remark}
\newtheorem{rem}[thm]{Remark}

\makeatletter
\let\c@equation\c@thm
\makeatother
\numberwithin{equation}{section}

\def\bthm{\begin{thm}}
\def\ethm{\end{thm}}
\def\blm{\begin{lem}}
\def\elm{\end{lem}}
\def\bdf{\begin{defn}}
\def\edf{\end{defn}}
\def\bpf{\begin{proof}}
\def\epf{\end{proof}}
\def\bpp{\begin{prop}}
\def\epp{\end{prop}}
\def\bcor{\begin{cor}}
\def\ecor{\end{cor}}
\def\brm{\begin{rem}}
\def\erm{\end{rem}}
\def\beg{\begin{exmp}}
\def\eeg{\end{exmp}}
\def\bD{\mathbb{D}}

\def\bG{\mathbb{G}}

\def\bN{\mathbb{N}}

\def\bQ{\mathbb{Q}}

\def\bZ{\mathbb{Z}}
\def\cA{\mathcal{A}}
\def\cB{\mathcal{B}}
\def\cC{\mathcal{C}}
\def\cD{\mathcal{D}}
\def\cE{\mathcal{E}}
\def\cF{\mathcal{F}}
\def\cG{\mathcal{G}}
\def\cH{\mathcal{H}}

\def\cK{\mathcal{K}}
\def\cL{\mathcal{L}}
\def\cM{\mathcal{M}}
\def\cN{\mathcal{N}}
\def\cO{\mathcal{O}}
\def\cP{\mathcal{P}}

\def\cR{\mathcal{R}}
\def\cS{\mathcal{S}}
\def\cT{\mathcal{T}}

\def\cX{\mathcal{X}}
\def\cY{\mathcal{Y}}
\def\cZ{\mathcal{Z}}

\def\scI{\mathscr{I}}

\newcommand{\raq}{\,\rightarrow \,}
\newcommand{\laq}{\,\leftarrow \,}

\newcommand{\xraq}[2][]{\, \xrightarrow[#1]{#2} \,}
\newcommand{\xlaq}[2][]{\, \xleftarrow[#1]{#2} \,}

\newcommand{\ra}{\rightarrow}
\newcommand{\la}{\leftarrow}
\newcommand{\rinto}{\hookrightarrow}

\newcommand{\ronto}{\twoheadrightarrow}

\newcommand{\rsa}{\stackrel{\sim}{\rightarrow}}

\newcommand{\xra}[2][]{\xrightarrow[#1]{#2}}
\newcommand{\xla}[2][]{\xleftarrow[#1]{#2}}

\newcommand{\ie}{{\it i.e.}}
\newcommand{\eg}{{\it e.g.}}

\newcommand{\Mod}{{\rm Mod}}

\newcommand{\Ab}{{\rm Ab}}

\newcommand{\Ch}{{\rm Ch}}
\newcommand{\cone}{{\rm cone}}
\newcommand{\cocone}{{\rm cocone}}

\newcommand{\Tor}{{\rm Tor}}
\newcommand{\Ext}{{\rm Ext}}

\newcommand{\Extcom}{\underline{{\rm Ext}}}

\newcommand{\Ob}{{\rm Ob}}
\newcommand{\op}{{\rm op}}

\newcommand{\id}{{\rm id}}

\newcommand{\Hom}{{\rm Hom}}
\newcommand{\Homcom}{\underline{{\rm Hom}}}

\newcommand{\coker}{{\rm coker}}
\newcommand{\colim}{{\rm colim}}

\newcommand{\cHom}{\mathscr{H}\text{\kern -3pt {\calligra\large om}}\,}

\newcommand{\RHom}{{\bm R}{\rm Hom}}
\newcommand{\RHomcom}{{\bm R} \underline{{\rm Hom}}}

\newcommand{\RcHom}{{\bm R}\cHom}
\newcommand{\RcHomcom}{{\bm R}\underline{\cHom}}
\newcommand{\cHomcom}{\underline{\cHom}}

\newcommand{\holim}{{\rm holim}}

\newcommand{\hocolim}{{\rm hocolim}}

\newcommand{\Spec}{{\rm Spec}}

\newcommand{\Proj}{{\rm Proj}}
\newcommand{\PProj}{\mathpzc{Proj}}

\newcommand{\QCoh}{{\rm QCoh}}
\newcommand{\Coh}{{\rm Coh}}

\newcommand{\Dsuit}{\cD^{\tiny \mbox{$\spadesuit $}}}

\newcommand{\Dbcoh}{\cD^b_{{\rm coh}}}

\newcommand{\Dqcoh}{\cD_{{\rm qcoh}}}

\newcommand{\Dmcoh}{\cD^{-}_{{\rm coh}}}

\newcommand{\coh}{{\rm coh}}

\newcommand{\Gr}{{\rm Gr}}
\newcommand{\GrA}{{\rm Gr}(A)}
\newcommand{\gr}{{\rm gr}}

\newcommand{\QpGr}{{\rm Q}^+{\rm Gr}}

\newcommand{\qpgr}{{\rm q}^+{\rm gr}}

\newcommand{\QpmGr}{{\rm Q}^+_{(m)}{\rm Gr}}
\newcommand{\QpdGr}{{\rm Q}^+_{(d)}{\rm Gr}}

\newcommand{\Torp}{{\rm Tor}^+}

\newcommand{\IIT}{I^{\infty}\text{-}{\rm Tor}}

\newcommand{\cIIT}{\scI^{\infty}\text{-}{\rm Tor}}

\newcommand{\Icomp}{I\text{-}{\rm comp}}

\newcommand{\RGam}{{\bm R}\Gamma}

\newcommand{\RGI}{{\bm R}\Gamma_{I}}

\newcommand{\tRGam}{\widetilde{{\bm R}\Gamma}}

\newcommand{\Ce}{\check{\cC}}
\newcommand{\CI}{\check{\cC}_I}

\newcommand{\LLam}{{\bm L}\Lambda}

\newcommand{\LLI}{{\bm L}\Lambda_{I}}

\newcommand{\ITR}{I\text{-}{\rm triv}}
\newcommand{\IpTR}{I^+\text{-}{\rm triv}}

\newcommand{\cITR}{\scI\text{-}{\rm triv}}

\newcommand{\GrMod}{{\rm GrMod}}

\newcommand{\lc}{{\rm lc}}

\DeclareMathAlphabet{\mathpzc}{OT1}{pzc}{m}{it}

\title{Grothendieck duality and Greenlees-May duality on graded rings}
\author{Wai-Kit Yeung}
\address{Department of Mathematics,
	Indiana University, Bloomington, IN 47405, USA}
\email{yeungw@iu.edu}
\begin{document}

\begin{abstract}
We formulate and prove Serre's equivalence for $\bZ$-graded rings. When restricted to the usual case of $\bN$-graded rings, our version of Serre's equivalence also sharpens the usual one by replacing the condition that $A$ be generated by $A_1$ over $A_0$ by a more natural condition, which we call the Cartier condition. For $\bZ$-graded rings coming from flips and flops, this Cartier condition relates more naturally to the geometry of the flip/flop in question. We also interpret Grothendieck duality as an instance of Greenlees-May duality for graded rings. These form the basic setting for a homological study of flips and flops in \cite{Yeu20a, Yeu20b}.
\end{abstract}

\maketitle


\tableofcontents

\section{Introduction}

One of the cornerstones in the study of projective varieties is that they bear a close relation with graded rings. 
In the usual formulation of this relation, one takes an $\bN$-graded ring $A$ such that $A_0$ is Noetherian, and $A$ is finitely generated by $A_1$ over $A_0$. Then there is a projective morphism $\pi : \Proj \, A \ra \Spec \, A_0$. 
%
Moreover, Serre's equivalence describes the category of coherent sheaves on $X = \Proj \, A$ as the Serre quotient $\Coh(X) \simeq \gr(A) / {\rm tor}(A)$ of the category $\gr(A)$ of finitely generated graded modules, by the Serre subcategory of torsion modules.
This allows for a study of the geometry of projective varieties via commutative algebra of graded rings.
In fact, since the Serre quotient $\gr(A) / {\rm tor}(A)$ is also well-defined for noncommutative graded rings, this opens up the possibilty of a study of ``noncommutative projective spaces'' (see, \eg, \cite{AZ94}).

Most of these studies of graded rings (commutative or not) have focused on the case of $\bN$-graded rings. 
There are other contexts, however, where the graded ring may not be concentrated in non-negative weights (\ie, these are $\bZ$-graded rings). One such context is from birational geometry, where there are two important classes of birational operations -- flips and flops -- both of which are controlled by (sheaves of) $\bZ$-graded rings.
Namely, if we are given a log flip%
\footnote{Roughly speaking, a log flip is a configuration \eqref{log_flip_intro} in which the $\pi^+$-ample direction is in a suitable sense aligned with the $\pi^-$-anti-ample direction. See, \eg, \cite{KM98} for details ({\it cf}. \cite{Tha96, Yeu20a})}
\begin{equation}  \label{log_flip_intro}
\begin{tikzcd} [row sep = 0]
X^-  \ar[rd,"\pi^-"] & &  X^+ \ar[ld, "\pi^+"'] \\
& Y
\end{tikzcd}
\end{equation}
then there is a sheaf $\cA$ of $\bZ$-graded rings on $Y$ such that $X^+ = \Proj^+_Y(\cA) := \Proj_Y(\cA_{\geq 0})$ and $X^- = \Proj^-_Y(\cA) := \Proj_Y(\cA_{\leq 0})$.


The present paper is devoted to establishing some basic results relating sheaves on $X^{+} := \Proj^{+}(A)$ with graded modules over $A$. By symmetry, these results also apply to $X^-$, and will form the foundations for a homological study of flips and flops in \cite{Yeu20a, Yeu20b}.

A convenient setting to express these relations is provided by Greenlees-May duality. By that we mean a certain recollement of derived categories
\begin{equation}  \label{GM_recollement_intro}
\begin{tikzcd} [row sep = 0, column sep = 50]
& & \cD_{\IIT}(\GrA) \ar[ld, bend right, "\iota"'] \ar[dd, shift left = 3, "\LLI"]  \ar[dd, phantom, "\simeq"] \\
\cD_{\ITR}(\GrA) \ar[r, "\iota"] 
& \cD(\GrA) \ar[ru, "\RGI"] \ar[rd, "\LLI"']  \ar[l, bend right, "\Ce_I"'] \ar[l, bend left, "\cE_I"]  \\
& & \cD_{\Icomp}(\GrA)  \ar[lu, bend left, "\iota"] \ar[uu, shift left = 3 , "\RGI"]
\end{tikzcd}
\end{equation}
associated to a graded ideal $I \subset A$. This is a direct extension of the usual Greenlees-May duality to graded rings. See Section \ref{GM_gr_rings_sec} for details.

When we take $I$ to be the graded ideal $I = I^+ := A_{>0} \cdot A$, this recollement is closely related to the projective space $X^+ = \Proj^+(A) := \Proj(A_{\geq 0})$.
Indeed, if we denote by $\pi^+ : \Proj^+(A) \ra \Spec(A_0) =: Y$ the canonical map, then there are three aspects of this relation:
\begin{equation}  \label{GGM_three_points}
\parbox{40em}{(1) The category $\cD_{\IpTR}(\GrA)$ is closely related to the derived category $\cD(\QCoh(X^+))$.\\
	(2) The functor $\Ce_{I^+}$ is closely related to the derived pushforward functor ${\bm R}\pi^+_*$.\\
	(3) The functor $\cE_{I^+}$ is closely related to the shriek pullback functor $(\pi^+)^!$.}
\end{equation}

The relation \eqref{GGM_three_points}(1) has two aspects. At the level of abelian categories, it is about Serre's equivalence. We sharpen the usual statement of Serre's equivalence (see Theorem \ref{Serre_equiv_thm}) so that the usual condition of $A_{\geq 0}$ being generated by $A_1$ over $A_0$ is replaced by a certain Cartier condition (Definition \ref{frac_Cartier}). In examples of $\bZ$-graded rings coming from log flips, this corresponds to the condition that the $\bQ$-Cartier divisors that define the log flip are Cartier (see \cite[Proposition 3.5]{Yeu20a}).

The second aspect of \eqref{GGM_three_points}(1) is the passage from abelian categories to derived categories. More precisely, Serre's equivalence identifies $\cD(\QCoh(X^+))$ as the derived category of a quotient, while the recollement \eqref{GM_recollement_intro} identifies $\cD_{\IpTR}(\GrA)$ as the quotient of derived categories. In Appendix \ref{Dcat_app}, we give some sufficient conditions that guarantee that the derived category of a Serre quotient coincides with the Verdier quotient of the corresponding derived categories. These conditions are satisfied in the present case, and we obtain the relation \eqref{GGM_three_points}(1)
(see the end of Section \ref{Serre_equiv_sec}).

The relation \eqref{GGM_three_points}(2) is the simplest: it essentially says that one can compute cohomology of quasi-coherent sheaves by the \v{C}ech complex (see \eqref{Ce_i_RGam}). 
The relation \eqref{GGM_three_points}(3) expresses Grothendieck duality in terms of Greenlees-May duality (see Theorem \ref{GGM_thm}).  
Both of these relations are crucial inputs in the homological study of flips and flops in \cite{Yeu20a}.

\section{Greenlees-May duality for graded rings}   \label{GM_gr_rings_sec}

In this section, we extend some results on the derived category of a commutative ring, especially Greenlees-May duality, to the graded case.

\bdf
A \emph{$\bZ$-graded ring} is a commutative ring $A$ with a $\bZ$-grading $A = \bigoplus_{n \in \bZ} A_n$. 
Here, by commutative we mean $xy = yx$, not $xy = (-1)^{|x||y|} yx$. 

A \emph{graded module} over $A$ will always mean a $\bZ$-graded module $M = \bigoplus_{n \in \bZ} M_n$.
\edf

We first recall the following result (see, e.g., \cite[Theorem 1.5.5]{BH93}):
\bpp  \label{Noeth_gr_ring}
Let $A$ be a $\bZ$-graded ring. The the followings are equivalent:
\begin{enumerate}
	\item $A$ is a Noetherian ring; 
	\item every graded ideal of $A$ is finitely generated;
	\item $A_0$ is Noetherian, and both $A_{\geq 0}$ and $A_{\leq 0}$ are finitely generated over $A_0$;
	\item $A_0$ is Noetherian, and $A$ is finitely generated over $A_0$.
\end{enumerate}
\epp

Denote by $\Gr(A)$ the category of graded modules over $A$, whose morphisms are maps of graded modules of degree $0$. 
Given two graded modules $M,N \in \Gr(A)$, then the $A$-module $M \otimes_A N$ has a natural grading where $\deg(x\otimes y) = \deg(x) + \deg(y)$ for homogeneous $x,y \in A$. Moreover, one can define a graded $A$-module $\Homcom_A(M,N)$ whose degree $i$ part is the set of $A$-linear homomorphism from $M$ to $N$ of homogeneous degree $i$. 
Thus, in particular, we have $\Hom_A(M,N) := \Hom_{\Gr(A)}(M,N) = \Homcom_A(M,N)_0$. These form the internal Hom objects with respect to the graded tensor product. 

The abelian category $\GrA$ is a Grothendieck category, with a set $\{A(-i)\}_{i \in \bZ}$ of generators. The same set is also a set of compact generators in the derived category $\cD(\GrA)$.
Since $\GrA$ is a Grothendieck category, the category of complexes has enough K-injectives (see, e.g., \cite[Tag 079P]{Sta}). Moreover, as in the ungraded case, it also has enough K-projectives (see, e.g., \cite[Tag 06XX]{Sta}).
As a result, the bifunctors $-\otimes_A - $ and $\Homcom_A(-,-)$ admit derived functors
\begin{equation*}
\begin{split}
- \otimes_{A}^{{\bm L}} - \, &: \, \cD(\GrA) \, \times \, \cD( \GrA) \raq \cD(\GrA) \\
\RHomcom_{A}(-,-) \, &: \, \cD(\GrA)^{\op} \, \times \, \cD( \GrA ) \raq \cD(\GrA)
\end{split}
\end{equation*}
which can in turn be used to define $\Extcom^{\bullet}_A(M,N)$ and $\Tor^A_{\bullet}(M,N)$.



\bdf  \label{I_infty_torsion_def}
Let $I$ be a graded ideal in a $\bZ$-graded ring $A$. Given any graded module $M$ over $A$, an element $x\in M$ is said to be \emph{$I^{\infty}$-torsion} if for every $f \in I$ there exists some $n > 0$ such that $f^n x = 0$. If $I$ is finitely generated, this is equivalent to $I^n x = 0$ for some $n > 0$. The graded module $M$ is said to be \emph{$I^{\infty}$-torsion} if every element in it is $I^{\infty}$-torsion.
Denote by $\IIT \subset \Gr(A)$ the full subcategory consisting of $I^{\infty}$-torsion modules.
\edf

It is clear that $\IIT \subset \GrA$ is a Serre subcategory. Thus the full subcategory $\cD_{\IIT}(\GrA) \subset \cD(\GrA)$ is a triangulated subcategory.
If $I$ is finitely generated, then this inclusion has a right adjoint, which has a simple and useful description.
To this end, we recall the following


\bdf  \label{RGam_f_def}
Let $f_1,\ldots,f_r$ be homogeneous elements in $A$, of degrees $d_1,\ldots,d_r$ respectively. 
For any graded module $M \in \Gr(A)$, 
we define the \emph{local cohomology complex} (or \emph{extended \v{C}ech complex}) of $M$ with respect to the tuple $(f_1,\ldots,f_r)$
to be the cochain complex of graded modules
\begin{equation}  \label{RGam_f}
\RGam_{(f_1,\ldots,f_r)}(M) \, := \, \bigl[ \, M \xra{d^0} \prod_{1 \leq i_0 \leq r} M_{f_{i_0}}  \xra{d^1}  
\prod_{1 \leq i_0 < i_1 \leq r} M_{f_{i_0}f_{i_1}}  \xra{d^2} \ldots \xra{d^{r-1}} 
M_{f_1\ldots f_r} \, \bigr]
\end{equation}
whose differentials are defined by 
$d^m := \sum_{j=0}^m (-1)^j d^m_j$, where $d^m_j$ is the direct product of the canonical maps $d^m_j : A_{i_0\ldots \hat{i_j} \ldots i_m} \ra A_{i_0\ldots i_m}$.
Here, the first term $M$ is put in cohomological degree $0$. 

For a cochain complex $M \in \Ch(\Gr(A))$ of graded modules, we define $\RGam_{(f_1,\ldots,f_r)}(M)$ to be the total complex of the double complex $C^{p,q} = \RGam_{(f_1,\ldots,f_r)}(M^p)^q$.
\edf

It is clear that we have
\begin{equation}   \label{RGam_tensor}
\RGam_{(f_1,\ldots,f_r)}(M) \, \cong \,  \RGam_{(f_1,\ldots,f_r)}(A) \otimes_A M \, \cong \,\RGam_{f_1}(A) \otimes_A \ldots \otimes_A \RGam_{f_r}(A) \otimes_A M
\end{equation}

The local cohomology complex may be written as a directed colimit of (cohomological) Koszul complexes
\begin{equation}  \label{Cech_Koszul}
\RGam_{(f_1,\ldots,f_r)}(M) \, \cong \, \colim_{(m_1,\ldots,m_r) \in (\bZ_{>0})^r} \, K^{\bullet}(M; f_1^{m_1},\ldots, f_r^{m_r})
\end{equation}
which can be computed, as is often done in the literature (see, e.g., \cite[Tag 0913]{Sta}, \cite[Theorem 3.5.6]{BH93} in the ungraded case), via the colimit on the cofinal system $(m,\ldots,m) \in (\bZ_{> 0})^r$.

Whenever a homogeneous element $g \in A$, say of degree $m$, lies in the ideal generated by $f_1^{m_1},\ldots, f_r^{m_r}$, the map in $\Ch(\Gr(A))$
\begin{equation*}
\cdot g \, : \, K^{\bullet}(M; f_1^{m_1},\ldots, f_r^{m_r}) \ra K^{\bullet}(M; f_1^{m_1},\ldots, f_r^{m_r})(m)
\end{equation*} 
is homotopic to zero, and hence induces the zero map in cohomology.
Indeed, the graded analogue of, say, \cite[Tag 0626]{Sta} establishes this for $M = A$, which then implies that it holds for all $M \in \Ch(\Gr(A))$, since we have 
$K^{\bullet}(M,f_1^{m_1},\ldots,f_r^{m_r}) = M \otimes_A K^{\bullet}(A,f_1^{m_1},\ldots,f_r^{m_r})$.
Since directed colimits commute with taking cohomology, we see by \eqref{Cech_Koszul} that the cohomology modules of the local cohomology complex are $I^{\infty}$-torsion, where $I = (f_1,\ldots,f_r)$ is the ideal generated by the elements $f_i$.

Since the complex $\RGam_{(f_1,\ldots,f_r)}(A)$ is flat over $A$, the functor 
$M \mapsto \RGam_{(f_1,\ldots,f_r)}(M)$ on $\Ch(\Gr(A))$ is exact, and hence descend to a functor at the level of derived categories. Moreover, we have seen that this functor has image inside the full subcategory $\cD_{\IIT}(\Gr (A))$. Thus, this gives a functor
\begin{equation}  \label{RGI_def}
\RGI := \RGam_{(f_1,\ldots,f_r)} \, : \, \cD(\Gr(A)) \raq \cD_{\IIT}(\Gr (A))
\end{equation}

Moreover, the map $\epsilon_M : \RGam_{(f_1,\ldots,f_r)}(M) \ra M$ defined by projecting to the first component of \eqref{RGam_f} gives rise to a natural transformation
\begin{equation}  \label{RGam_counit}
\epsilon \, : \, \iota \circ \RGI \, \Rightarrow \, \id
\end{equation}  
where $\iota : \cD_{\IIT}(\Gr (A)) \ra \cD(\Gr(A))$ is the inclusion functor.
The cone of $\epsilon_M$ is homotopic to the kernel of $\epsilon_M$ shifted by $1$, which is given by the following
\bdf  \label{Cech_cplx_def}
The \emph{\v{C}ech complex} of a graded module $M$ with respect to a tuple $(f_1,\ldots,f_r)$ of homogeneous elements
is the cochain complex of graded modules
\begin{equation}  \label{Ce_f}
\Ce_{(f_1,\ldots,f_r)}(M) \, := \, \bigl[ \, \prod_{1 \leq i_0 \leq r} M_{f_{i_0}}  \xraq{-d^1}  
\prod_{1 \leq i_0 < i_1 \leq r} M_{f_{i_0}f_{i_1}}  \xraq{-d^2} \ldots \xraq{-d^{r-1}} 
M_{f_1\ldots f_r} \, \bigr]
\end{equation}
given as a subcomplex of \eqref{RGam_f}, 
shifted by one.
As in Definition \ref{RGam_f_def}, this definition can be extended to cochain complexes $M \in \Ch(\Gr(A))$ by taking the total complex.
\edf

Clearly, the natural transformation \eqref{RGam_counit} is an isomorphism for $M \in \cD(\Gr(A))$ if and only if $\Ce_{(f_1,\ldots,f_r)}(M)$ has zero cohomology. 
Now if $M \in \cD_{\IIT}(\Gr(A))$, then each of the terms in \eqref{Ce_f}, thought of as a column in a double complex, has zero cohomology. Thus, as an iterative cone of complexes with zero cohomology, 
the total complex $\Ce_{(f_1,\ldots,f_r)}(M)$ also has zero cohomology. This shows that
\begin{equation}  \label{tor_epsilon_isom}
\parbox{40em}{If $M \in \cD_{\IIT}(\Gr(A))$, then the natural transformation $\epsilon_M$ is an isomorphism in $\cD(\Gr(A))$.}
\end{equation}

As a formal consequence of \eqref{RGI_def}, \eqref{RGam_counit} and \eqref{tor_epsilon_isom}, we have the following
\bthm   \label{RGam_right_adj}
The functor \eqref{RGI_def} is a right adjoint to the inclusion $\iota : \cD_{\IIT}(\Gr (A)) \ra \cD(\Gr(A))$, with counit given by \eqref{RGam_counit}.
\ethm

For each $M \in \cD(\GrA)$, there is an exact triangle
\begin{equation}  \label{RGam_Ce_seq}
\ldots \raq \RGI(M) \xraq{\epsilon_M} M \xraq{\eta_M} \CI(M) \xraq{\delta_M} \RGI(M)[1] \raq \ldots 
\end{equation}
where $\eta_M = -d^0$, the negative of the first differential in \eqref{RGam_f},
and $\delta_M$ is the inclusion.

If $A$ is Noetherian, then the functor $\RGI$ has an alternative descriptions as right derived functor.
Namely, let $\Gamma_I : \Gr(A) \ra \IIT$ be the right adjoint to the inclusion $\IIT \ra \Gr(A)$, then we have the following graded analogue of a well-known statement (see, e.g., \cite[Tag 0955]{Sta}), which justifies the name ``local cohomology complex'' of $\RGam_I(M)$:

\bpp  \label{local_cohom_der_func}
If $A$ is Noetherian, then the canonical functor $\cD(\IIT) \ra \cD_{\IIT}(\Gr(A))$ is an equivalence. Moreover, under this equivalence, the functor \eqref{RGI_def} is identified with the right derived functor of $\Gamma_I : \Gr(A) \ra \IIT$.
\epp

As a consequence, the functor $\Ce_I$ also has an interpretation as a right derived functor. 
Namely, let $\!\,^0 \Ce_I : \GrA \ra \GrA$ be the functor $\!\,^0 \Ce_I(M) := H^0(\Ce_I(M))$, then we have

\bcor  \label{Cech_der_func}
If $A$ is Noetherian, then the functor $\Ce_I : \cD(\GrA) \ra \cD(\GrA)$ is the right derived functor of  $\,^0 \Ce_I : \GrA \ra \GrA$.
\ecor

\bpf
Notice that for any $M \in \GrA$, there is an exact sequence 
\begin{equation}  \label{Gam_M_Ce_four_terms}
0 \raq \Gamma_I(M) \raq M \raq \!\,^0 \Ce_I(M) \raq H^1( \RGI(M) ) \raq 0
\end{equation}
Now if we apply this termwise to a K-injective complex of injective graded modules $M^{\bullet}$, then Proposition \ref{local_cohom_der_func} shows that $H^1( \RGI(M^{i}) ) = 0$ for all $i \in \bZ$, and we therefore have a short exact sequence of cochain complexes 
\begin{equation*}
0 \raq \Gamma_I(M^{\bullet}) \raq M^{\bullet} \raq \!\,^0 \Ce_I(M^{\bullet}) \raq 0
\end{equation*}
By Proposition \ref{local_cohom_der_func} again, the map $\Gamma_I(M^{\bullet}) \raq M^{\bullet}$ here is isomorphic in $\cD(\GrA)$ to the map $\epsilon_M : \RGI(M^{\bullet}) \ra M^{\bullet}$ in \eqref{RGam_Ce_seq}. Thus, the cone of these two maps are also identified. In other words, we have $\!\,^0 \Ce_I(M^{\bullet}) \cong \Ce_I(M^{\bullet})$ in $\cD(\GrA)$.
\epf


%
%

\vspace{0.3cm}

Now we discuss a closely related notion of derived complete graded modules, following \cite{Sta}.
Let $M \in \Ch(\Gr(A))$ be a chain complex, and let $f \in A$ be a homogeneous element, say of degree $m$. Denote by $T(M,f)$ the cochain complex in $\Gr(A)$ defined as the homotopy limit (see \eqref{holim_def})
\begin{equation}  \label{TMf_def}
T(M,f) \, := \, \holim \, [ \, M \xlaq{f} M(-m) \xlaq{f} M(-2m) \xlaq{f} \ldots  \, ]
\end{equation}

This cochain complex can be written as a Hom-complex. Indeed, if we denote by $A_{\{f\}}$ the homotopy colimit (see \eqref{hocolim_def})
\begin{equation}  \label{A_brac_f}
A_{\{f\}} \, := \, \hocolim \, [ \, A \xraq{f} A(m) \xraq{f} A(2m) \xraq{f} \ldots  \, ]
\end{equation}
then, since the corresponding ordinary colimit is simply $A_f$, we have, by Lemma \ref{hocolim_colim_qism}, a quasi-isomorphism $A_{\{f\}} \xra{\simeq} A_f$,
thus giving a free resolution $A_{\{f\}}$ of $A_f$.
This allows us to rewrite \eqref{TMf_def} as
\begin{equation}  \label{TMf_RHom}
T(M,f) \, = \, \Homcom_A( A_{\{f\}} , M ) \, \simeq \, \RHomcom_A( A_f , M )
\end{equation}

Since the association $M \mapsto T(M,f)$ is exact, it descends to a functor $T(-,f) : \cD(\Gr(A)) \ra \cD(\Gr(A))$.
We start with the following graded analogue of \cite[Tag 091P]{Sta}:
\blm  \label{TMf_equiv}
Let $f \in A$ be a homogeneous element of degree $d$. 
Then for any complex $M \in \cD(\Gr(A))$, the followings are equivalent:
\begin{enumerate}
	\item $T(M,f)$ has zero cohomology;
	\item $\RHomcom_A(E,M)$ has zero cohomology for all $E \in \cD(\Gr(A_f))$;
	\item for every $p \in \bZ$, we have $\Homcom_A(A_f, H^p(M)) = 0$ and $\Extcom^1_A(A_f, H^p(M)) = 0$.
	\item for every $p \in \bZ$, the complex $T(H^p(M),f)$ has zero cohomology.
\end{enumerate}
\elm

\bpf
Since $\{ A_f(n)\}_{n \in \bZ}$ is a set of compact generators of $\cD(\Gr(A_f))$, we see that $\RHomcom(E,M) \simeq 0$ for all $E \in \cD(\Gr(A_f))$ if and only if it holds for $E = A_f$. In view of \eqref{TMf_RHom}, this proves the equivalence $(1) \Leftrightarrow (2)$.
The equivalence $(3) \Leftrightarrow (4)$ is also obvious, since the graded modules appearing in (3) are simply the cohomology of the complex $T(H^p(M),f)$ in (4).
For the equivalence $(1)\Leftrightarrow (3)$, simply take a spectral sequence.
\epf

\brm  \label{mult_f_isom_rem}
Since $A_f$ is flat over $A$, the standard tensor-forgetful adjunction between $\cD(\Gr(A_f))$ and $\cD(\Gr(A))$ identifies the former as the full subcategory of $\cD(\Gr(A))$ consisting of complexes such that the multiplication map $f : E \ra E(n)$ is an isomorphism in $\cD(\GrA)$. Thus condition (2) of Lemma \ref{TMf_equiv} may be rewritten with respect to such complexes $E$.
\erm

Given $M \in \cD(\Gr(A))$, let $I_M$ be the subset of $A$ 
consisting of elements $f = \sum f_i \in A$ such that, for each homogeneous component $f_m \in A_m$ of $f$, 
the complex $T(M,f_m)$ has zero cohomology. Then exactly the same proof as in \cite[Tag 091Q]{Sta} shows the following
\blm
The subset $I_M \subset A$ is a radical graded ideal.
\elm

\bdf  \label{derived_comp_def}
Let $I \subset A$ be a graded ideal. An object $M \in \cD(A)$ is said to be \emph{derived complete with respect to $I$} if for every homogeneous $f \in I$, the complex $T(M,f)$ has zero cohomology. In other words, if $I \subset I_M$. 

A graded module $M \in \Gr(A)$ is said to be \emph{derived complete with respect to $I$} if the corresponding complex $M[0] \in \cD(\Gr(A))$ concentrated in cohomological degree zero is so.
\edf

Denote by $\cD_{\Icomp}(\GrA) \subset \cD(\GrA)$ the full subcategory consisting of objects that are derived complete with respect to $I$.
Since $T(M,f)$ can be written as a derived Hom complex \eqref{TMf_RHom}, this full subcategory 
is triangulated and split-closed. 
Moreover, by Lemma \ref{TMf_equiv}(4), an object $M \in \cD(\GrA)$ is in $\cD_{\Icomp}(\GrA)$
if and only if all of its cohomology modules are in the subcategory $\Icomp \subset \GrA$. 
This, in turn, shows that $\Icomp$ is a weak Serre subcategory.

Now we focus on the case when $I$ is finitely generated, say by $f_1,\ldots,f_r$. 
In this case, we first construct a free resolution of the local cohomology complex $\RGam_{(f_1,\ldots,f_r)}(A)$ by replacing each term $A_{f_{i_0}\ldots f_{i_p}}$ in \eqref{RGam_f}
by its free resolution $A_{\{ f_{i_0}\ldots f_{i_p} \}}$ given in \eqref{A_brac_f}.
Indeed, for any two homogeneous elements $f$ and $g$, say of degrees $m$ and $n$ respectively, the canonical map $A_f \ra A_{fg}$ lifts canonically to a map $A_{\{f\}} \ra A_{\{fg\}}$ of cochain complexes, induced by a map of the corresponding directed systems defining the respective homotopy colimit, where the component $A(mi)$ in \eqref{A_brac_f} is sent to the component $A((m+n)i)$ via the map $g^i$.
Thus we define
\begin{equation}  \label{tRGam_f}
\tRGam_{(f_1,\ldots,f_r)}(A) \, := \, \bigl[ \, A \xra{\widetilde{d^0}} \prod_{1 \leq i_0 \leq r} A_{\{f_{i_0}\}}  \xra{\widetilde{d^1}} 
\prod_{1 \leq i_0 < i_1 \leq r} A_{\{f_{i_0}f_{i_1}\}}  \xra{\widetilde{d^2}} \ldots \xra{\widetilde{d^{r-1}}} 
A_{\{f_1\ldots f_r\}} \, \bigr]
\end{equation}
understood as a total complex of a double complex. This gives a free resolution
\begin{equation}  \label{tRGam_RGam_qism}
\tRGam_{(f_1,\ldots,f_r)}(A) \xraq{\simeq} \RGam_{(f_1,\ldots,f_r)}(A)
\end{equation}


\bdf
For any cochain complex $M \in \Ch(\GrA)$, define the cochain complex $\LLam_{(f_1,\ldots,f_r)}(M) \in \Ch(\GrA)$ to be the Hom-complex
\begin{equation*}
\begin{split}
&\LLam_{(f_1,\ldots,f_r)}(M) \, := \, \Homcom_A( \, \tRGam_{(f_1,\ldots,f_r)}(A) \, , \, M \, ) \\
& = \, \bigl[ \, M \xlaq{d_0} \prod_{1 \leq i_0 \leq r} T(M,f_{i_0})  \xlaq{d_1}  
\prod_{1 \leq i_0 < i_1 \leq r} T(M,f_{i_0}f_{i_1})  \xlaq{d_2} \ldots \xlaq{d_{r-1}} 
T(M,f_1\ldots f_r) \, \bigr]
\end{split}
\end{equation*}
where the last line is thought of as the total complex of a double complex.
\edf


In view of \eqref{tRGam_RGam_qism}, 
it gives an explicit model for the derived Hom complex
\begin{equation*}
\LLam_{(f_1,\ldots,f_r)}(M) \, \simeq \, \RHomcom_A( \RGam_{I}(A) , M )
\end{equation*}
where $I = (f_1,\ldots,f_r)$ is the ideal genereated by the elements $f_i$.
Moreover, recall that $\RGam_{I}(A)$ has $I^{\infty}$-torsion cohomology, so that $A_f \otimes_A \RGam_{I}(A) \simeq 0$ for all homogeneous elements $f \in I$.
As a result, we have
\begin{equation*}
T( \LLam_{(f_1,\ldots,f_r)}(M) , f ) \, \simeq \, \RHomcom_A( A_f , \RHomcom_A( \RGam_{I}(A) , M ) )\, \simeq \, \RHomcom_A( A_f \otimes_A \RGam_{I}(A) , M ) \simeq 0
\end{equation*}
for any homogeneous element $f \in I$. In other words, the complex $\LLam_{(f_1,\ldots,f_r)}(M)$ is always derived $I$-complete. Thus, the exact functor $\LLam_{(f_1,\ldots,f_r)}$ descends to a functor
\begin{equation}  \label{LLI_fucntor}
\LLI \, := \, \LLam_{(f_1,\ldots,f_r)} \, : \, \cD(\GrA) \raq \cD_{\Icomp}(\GrA)
\end{equation}

Moreover, the canonical map $\tRGam_{(f_1,\ldots,f_r)}(A) \ra A$ defined by projecting to the first component in \eqref{tRGam_f} induces a map $\epsilon^{\vee}_M : M \ra \LLam_{(f_1,\ldots,f_r)}(M)$ of cochain complexes. Thus, there is a natural transformation
\begin{equation}  \label{LLI_unit}
\epsilon^{\vee} \, : \, \id \, \Rightarrow \, \iota \circ \LLI
\end{equation}
where $\iota : \cD_{\Icomp}(\GrA) \ra \cD(\GrA)$ is the inclusion functor.
The cone of $\epsilon^{\vee}_M$ is homotopic to the cokernel of the (injective) map $\epsilon^{\vee}_M$. If we shift this cokernel by $-1$, we obtain the following

\bdf
Given a cochain complex $M \in \GrA$, define 
\begin{equation}  \label{cE_f}
\cE_{(f_1,\ldots,f_r)}(M)  :=  \bigl[ \, \prod_{1 \leq i_0 \leq r} T(M,f_{i_0})  \xla{-d^1}  
\prod_{1 \leq i_0 < i_1 \leq r} T(M,f_{i_0}f_{i_1})  \xla{-d^2} \ldots \xla{-d^{r-1}} 
T(M,f_1\ldots f_r) \, \bigr]
\end{equation}
understood as the total complex of a double complex.
\edf

Clearly the natural transformation \eqref{LLI_unit} is an isomorphism for $M \in \cD(\GrA)$ if and only if $\cE_{(f_1,\ldots,f_r)}(M)$ has zero cohomology. 
Now if $M \in \cD_{\Icomp}(\GrA)$ then each of the terms in \eqref{cE_f}, thought of as a column in a double complex, has zero cohomology. 
As an iterated cone on such columns, the total complex $\cE_{(f_1,\ldots,f_r)}(M)$ also has zero cohomology. This shows that
\begin{equation}  \label{eta_M_qism_comp}
\parbox{40em}{If $M \in \cD_{\Icomp}(\GrA)$, then the natural transformation $\epsilon^{\vee}_M$ is an isomorphism in $\cD(\GrA)$.}
\end{equation}

As a formal consequence of \eqref{LLI_fucntor}, \eqref{LLI_unit} and \eqref{eta_M_qism_comp}, we have the following

\bthm  \label{LLI_adj}
The functor \eqref{LLI_fucntor} is left adjoint to the inclusion $\iota : \cD_{\Icomp}(\GrA) \ra \cD(\GrA)$, with unit given by \eqref{LLI_unit}.
\ethm

\bdf
The complex $\LLam_I(M)$ is either called the \emph{local homology complex}, or the \emph{derived completion}, of $M$ with respect to the graded ideal $I \subset A$.
\edf

Dual to \eqref{RGam_Ce_seq}, we also have the following exact triangle
\begin{equation}  \label{LLI_cE_seq}
\ldots \raq \cE_I(M) \xraq{\eta_M^{\vee}} M \xraq{\epsilon^{\vee}_M} \LLI(M) \xraq{-\delta^{\vee}_M} \cE_I(M)[1] \raq \ldots 
\end{equation}

In fact, \eqref{RGam_Ce_seq} and \eqref{LLI_cE_seq} are dual in a very precise sense. Namely if we substitute $M = A$ in \eqref{RGam_Ce_seq} to get the exact triangle
\begin{equation}  \label{RGam_Ce_A_seq}
\ldots \raq \RGI(A) \xraq{\epsilon_A} A \xraq{\eta_A} \CI(A) \xraq{\delta_A} \RGI(A)[1] \raq  \ldots 
\end{equation}
then we have
\begin{equation*}
\eqref{RGam_Ce_seq} \, \simeq \, \eqref{RGam_Ce_A_seq} \otimes_A^{{\bm L}} M  \qquad \text{ and } \qquad 
\eqref{LLI_cE_seq} \, \simeq \, \RHomcom_A( \eqref{RGam_Ce_A_seq} , M  )
\end{equation*}

We record the following graded analogue of \cite[Tag 0A6S, 0A6D]{Sta}:

\blm  \label{mult_f_isom}
Given $E \in \cD(\GrA)$, suppose that there exists a homogeneous element $f \in I$, say of degree $n$, such that the multiplication map $f : E \ra E(n)$ is an isomorphism in $\cD(\GrA)$, then we have $\RGI(E) \simeq 0$ and $\LLI(E) \simeq 0$.
\elm

\bpf
The property that $f : E \ra E(n)$ is an isomorphism in $\cD(\GrA)$ is preserved under $\RGI(A)\otimes_A^{{\bm L}} - $, so that $\RGI(E)$ also has this property. However, cohomology modules of $\RGI(E)$ are $I^{\infty}$-torsion, so that in particular any element $x \in H^p(\RGI(E))$ is annihilated by some high enough powers of $f$. Hence we must have $x = 0$.

To show that $\LLI(E) \simeq 0$, notice that, by the criterion in Lemma \ref{TMf_equiv}(2) (see also Remark \ref{mult_f_isom_rem}), we have $\RHomcom_A(E,M) \simeq 0$ for all $M \in \cD_{\Icomp}(\GrA)$. Since the derived $I$-completion functor is left adjoint to inclusion (see Theorem \ref{LLI_adj}), this shows that $\LLI(E) \simeq 0$.
\epf

For any homogeneous element $f \in I$, the module $A_f$ clearly satisfies the condition of this Lemma. The same is therefore true for $M_f = A_f \otimes_A M$ and $T(M,f) \simeq \RHomcom_A(A_f,M)$. This observation allows us to prove the following
\bpp  \label{I_triv_equiv}
For any $M \in \cD(\GrA)$, we have $\LLI(\Ce_I(M)) \simeq 0$ and $\RGI(\cE_I(M)) \simeq 0$. Moreover,
\begin{equation*}
\RGI(M) \simeq 0 
\quad \Leftrightarrow \quad  \LLI(M) \simeq 0
\end{equation*}
\epp

\bpf
For the first statement, 
notice that in the double complexes \eqref{Ce_f} and \eqref{cE_f} defining $\Ce_I(M)$ and $\cE_I(M)$, each column satisfies the condition of Lemma \ref{mult_f_isom} (see the paragraph preceeding the present Proposition), so that we have $\LLI(\Ce_I(M)) \simeq 0$ and $\RGI(\cE_I(M)) \simeq 0$.

For the second statement, if $\RGI(M) \simeq 0 $, then we have $M \simeq \Ce_I(M)$ by \eqref{RGam_Ce_seq}, so that $\LLI(M) \simeq 0$ by the first statement. The converse is completely symmetric, using \eqref{LLI_cE_seq} this time.
\epf

\bcor
For any $M \in \cD(\GrA)$, we have the following isomorphisms in $\cD(\GrA)$:
\begin{equation*}
\RGI( M ) \xraq[\simeq]{ \RGI( \eta_M ) } \RGI( \LLI(M) ) 
\qquad \text{ and }\qquad 
\LLI( \RGI(M) ) \xraq[\simeq]{ \LLI( \epsilon_M ) } \LLI( M )
\end{equation*}
\ecor

\bdf  \label{ITR_def}
A complex $M \in \cD(\GrA)$ is said to be \emph{$I$-trivial} if we have $\RGI(M) \simeq 0$, or equivalently $\LLI(M) \simeq 0$ by Proposition \ref{I_triv_equiv}.
Denote by $\cD_{\ITR}(\GrA)$ the full subcategory of $\cD(\GrA)$ consisting of $I$-trivial complexes.
\edf

Thus, we have the following recollement (for two functors pointing in opposite horizontal directions, the functor on top of the other is implicitly understood to be the left adjoint of the other):
\begin{equation*}
\begin{tikzcd} [row sep = 0, column sep = 50]
& & \cD_{\IIT}(\GrA) \ar[ld, bend right, "\iota"'] \ar[dd, shift left = 3, "\LLI"]  \ar[dd, phantom, "\simeq"] \\
\cD_{\ITR}(\GrA) \ar[r, "\iota"] 
& \cD(\GrA) \ar[ru, "\RGI"] \ar[rd, "\LLI"']  \ar[l, bend right, "\Ce_I"'] \ar[l, bend left, "\cE_I"]  \\
& & \cD_{\Icomp}(\GrA)  \ar[lu, bend left, "\iota"] \ar[uu, shift left = 3 , "\RGI"]
\end{tikzcd}
\end{equation*}

In particular, this gives rise to two semi-orthogonal decompositions
\begin{equation}  \label{GM_two_SODs}
\begin{split}
\cD(\GrA) \, &= \, \langle \, \cD_{\ITR}(\GrA) \, , \,  \cD_{\IIT}(\GrA) \, \rangle \\
\cD(\GrA) \, &= \, \langle \, \cD_{\Icomp}(\GrA) \, , \, \cD_{\ITR}(\GrA)  \, \rangle
\end{split}
\end{equation}
whose decomposition triangles are \eqref{RGam_Ce_seq} and \eqref{LLI_cE_seq} respectively.

\vspace{0.3cm}

The notion of derived completeness is closely related to the usual notion of completeness with respect to a graded ideal. We will only mention one aspect of this relation, in the form of Proposition \ref{comp_der_comp_1} below.
Let $A$ be a $\bZ$-graded ring, and $I \subset A$ be a graded ideal. 
Then we may form the inverse system $(A/I^n)_{n \geq 1}$ of graded rings, and take its inverse limit
$\hat{A} := \varprojlim A/I^n$ in the category of graded rings,
which is computed as a limit in each graded component. 
Similarly, for any graded module $M$ over $A$, one can also take the 
inverse limit $\hat{M} := \varprojlim M/I^n M$ in the category of graded modules, which is again computed as a limit in each graded component. 
As usual, completions come with canonical maps $A \ra \hat{A}$ and $M \ra \hat{M}$.
We say that a graded module $M \in \Gr(A)$ is \emph{$I$-separated} if the canonical map $M \ra \hat{M}$ is injective; it is \emph{$I$-precomplete} if  $M \ra \hat{M}$ is surjective;
it is \emph{$I$-adically complete} if  $M \ra \hat{M}$ is bijective.
Then we have the following graded analogue of \cite[Tag 091R]{Sta}:


%
\bpp  \label{comp_der_comp_1}
Given a graded module $M \in \Gr(A)$, then we have
\begin{enumerate}
	\item If $M$ is $I$-adically complete, then it is derived complete with respect to $I$.
	\item If $I$ is finitely generated, and if $M$ is derived complete with respect to $I$,then $M$ is $I$-precomplete.
\end{enumerate} 
\epp

\section{Serre's equivalence}  \label{Serre_equiv_sec}

Let $A$ be a $\bZ$-graded ring. Denote by $I^+$ and $I^-$ the graded ideals 
$I^+ = A_{> 0}\cdot A$ and $I^- = A_{< 0}\cdot A$.
Let $X^+ = \Proj^+(A) := \Proj(A_{\geq 0})$. Thus, $X^+$ is covered by the standard open sets $D(f) \cong \Spec((A_{\geq 0})_{(f)})$, for homogeneous $f$ of positive degree. Since the canonical map $(A_{\geq 0})_{f} \ra A_{f}$ is an isomorphism of graded rings for $\deg(f) > 0$, the standard open sets $D(f)$ can also be described as $D(f) \cong \Spec(A_{(f)})$. 

In the study of the scheme $\Proj^+(A)$, we will mostly consider the case when $A_{\geq 0}$ is Noetherian. In view of Proposition \ref{Noeth_gr_ring}, we may make the following more precise assumption:
\begin{equation}  \label{I_plus_gen}
\parbox{40em}{Assume that $A_{\geq 0}$ is generated over $A_0$ by homogeneous elements $f_1,\ldots,f_p$ of positive degrees $d_i := \deg(f_i) > 0$, and let $d > 0$ be a positive integer that is divisible by each of $d_i$.}
\end{equation}
Then we record the following elementary result:
\blm  \label{AN_Ad}
Under the assumption \eqref{I_plus_gen},
for any $N > dp - \sum_{i=1}^p d_i$, we have $A_{N} = A_d \cdot A_{N-d}$.
\elm
\bpf
Given $f_1^{e_1}\ldots f_p^{e_p}g_1^{l_1}\ldots g_q^{l_q} \in A_{N}$ with $\deg(g_i) \leq 0$, we must have $e_i \geq d/d_i$ for some $1 \leq i \leq p$. Then $f_1^{e_1}\ldots f_p^{e_p}g_1^{l_1}\ldots g_q^{l_q} = (f_i^{d/d_i})(f_1^{e_1}\ldots f_i^{e_i - d/d_i} \ldots f_p^{e_p}g_1^{l_1}\ldots g_q^{l_q})$.
\epf

This Lemma can be applied to give an alternative description of $(I^+)^{\infty}$-torsion elements in the sense of Definition \ref{I_infty_torsion_def}:

\blm  \label{IIT_element_equiv}
Under the assumption \eqref{I_plus_gen},
suppose we are given a homogeneous element $x$ in a graded module $M$, then the followings are equivalent:
\begin{enumerate}
	\item $x$ is $(I^+)^{\infty}$-torsion;
	\item $x \cdot (A_d)^n = 0$ for all sufficiently large $n$.
	\item $x \cdot A_{\geq s} = 0$ for some $s \in \bZ$.
\end{enumerate}
\elm

\bpf
Since $A_d \subset I^+$, the implication $(1) \Rightarrow (2)$ is immediate. 
To show that $(2) \Rightarrow (3)$, notice that Lemma \ref{AN_Ad} implies that $A_N = (A_d)^n \cdot A_{N-dn}$ for $N \geq dp + dn$, so that we have $x \cdot A_{\geq dp+dn} = 0$.
The implication $(3) \Rightarrow (1)$ is obvious since the ideal $I^+$ is generated by elements of positive degrees.
\epf

Each graded module $M \in \GrA$ is in particular a graded module over $A_{\geq 0}$, and hence gives an associated quasi-coherent sheaf $\widetilde{M}$ on $\Proj^+(A)$. 
In fact, it is more natural to consider $\widetilde{M}$ as piecing together the modules $M_{(f)}$ over $A_{(f)} = (A_{\geq 0})_{(f)}$, so that we may avoid the unnatural procedure of passing to $A_{\geq 0}$.
In this way, we see that the association $M \mapsto \widetilde{M}$ is lax monoidal, meaning that there are canonical maps $\widetilde{M} \otimes_{\cO_{X^+}} \widetilde{N} \ra \widetilde{M\otimes_A N}$,
satisfying the usual associativity and unitality conditions. 
Moreover, there is a canonical map of $A_0$-modules
\begin{equation}  \label{M_to_global_tilde_M}
M_0 \ra H^0(X^+ , \widetilde{M})
\end{equation}

Let $\widetilde{A}$ be the sheaf of $\bZ$-graded algebras on $\Proj^+(A)$ given by $\widetilde{A}_i := \widetilde{A(i)}$. The graded algebra structure is induced by the lax monoidal structure $\widetilde{A(i)} \otimes_{\cO_{X^+}} \widetilde{A(j)} \ra \widetilde{A(i+j)}$. In general, these maps may not be isomorphisms. However, we have the following
\blm  \label{Proj_A_twists}
Assume that \eqref{I_plus_gen} holds, then we have the followings:
\begin{enumerate}
	\item The sheaf $\widetilde{A(d)}$ is an ample%
	\footnote{We use the definition of ampleness in \cite[Tag 01PS]{Sta}.} invertible sheaf on $\Proj^+(A)$.
	\item For any $M \in \Gr(A)$ and any $i \in \bZ$, the map $\widetilde{M} \otimes_{\cO_{X^+}} \widetilde{A(di)} \ra \widetilde{M(di)}$ is an isomorphism.
\end{enumerate}
\elm

\bpf
By Lemma \ref{AN_Ad}, the standard open sets $D(f)$, for $\deg(f) = d$, covers $\Proj^+(A)$. Thus, condition (2) can be checked on $D(f)$, which is obvious. This also shows that $\widetilde{A(d)}$ is invertible. Finally, it is ample again because the standard open sets $D(f)$, for $\deg(f) = d$, covers $\Proj^+(A)$.
\epf

In fact, the condition \eqref{I_plus_gen} is not necessary to guarantee the conclusions of Lemma \ref{Proj_A_twists}. Instead, we may consider the following notion:

\bdf  \label{frac_Cartier}
A $\bZ$-graded ring $A$ is said to be \emph{positively $\tfrac{1}{d}$-Cartier}, for an integer $d > 0$, if 
the canonical map $\widetilde{A(di)} \otimes_{\cO_{X^+}} \widetilde{A(dj)} \ra \widetilde{A(di+dj)}$ is an isomorphism for all $i,j \in \bZ$.
Similarly, it is said to be \emph{negatively $\tfrac{1}{d}$-Cartier} if the analogous condition holds for $X^- = \Proj^-(A)$ in place of $X^+$.
It is said to be \emph{$\tfrac{1}{d}$-Cartier} if it is both positively and negatively $\tfrac{1}{d}$-Cartier.
When $d=1$, we simply say that $A$ is \emph{(positively/negatievly) Cartier}.
\edf

\blm
If $A$ is positively $\tfrac{1}{d}$-Cartier, then 
\begin{enumerate}
	\item $\Proj^+(A)$ can be covered by the standard open subsets $D(fg)$, where $f , g$ are homogeneous elements of positive degrees such that  $\deg(f) - \deg(g) = d$.
	\item $\widetilde{A(di)}$ is an invertible sheaf for all $i \in \bZ$, and the canonical map $\widetilde{M} \otimes_{\cO_{X^+}} \widetilde{A(di)} \ra \widetilde{M(di)}$ is an isomorphism for any $M \in \GrA$.
\end{enumerate}
\elm

\bpf
For the first statement, we may replace $A$ by $A_{\geq 0}$, so that it follows from \cite[Tag 01MU]{Sta}. The second statement clearly follows from the first.
\epf

For each $m > 0$, let $A^{(m)}$ be the $\bZ$-graded ring $(A^{(m)})_i = A_{mi}$. Then there is a canonical isomorphism $\Proj^+(A^{(m)}) \cong \Proj^+(A)$ given by $(A^{(m)})_{(f^m)} = A_{(f)}$.
For any graded module $M \in \GrA$, one also has an associated graded module $M^{(m)} \in \Gr(A^{(m)})$, defined by $(M^{(m)})_i := M_{mi}$.
This often allows one to deduce results on all weight components $i \in \bZ$ from the corresponding result on weight components $i \in m\bZ$, for some $m > 0$. 
This is the method we use to establish Proposition \ref{X_isom_Proj_2} below. To this end, we first recall the following standard result (see, e.g., \cite[Tags 01Q1, 01QJ]{Sta}):

\bpp  \label{X_isom_Proj_1}
Let $\pi : X \ra Y = \Spec \, R$ be a proper morphism, where $R$ is Noetherian, and let $\cL$ be an ample invertible sheaf on $X$. Define the $\bZ$-graded algebra $S$ by $S_i := H^0(X;  \cL^{\otimes i})$, then $S_{\geq 0}$ is Noetherian, and there is a canonical isomorphism $\varphi : X \xra{\cong} \Proj^+(S)$ over $\Spec \, R$.

Moreover, for any quasi-coherent sheaf $\cF$ on $X$, let $M \in \Gr(S)$ be the graded $S$-module given by $M_i := H^0(X, \cF \otimes \cL^{\otimes i})$, then there is a canonical isomorphism $\varphi^*( \widetilde{M} ) \xra{\cong} \cF$, such that the composition 
\begin{equation*}
M_0 \xraq{\eqref{M_to_global_tilde_M}} H^0(\Proj^+(S), \widetilde{M}) \, = \,  H^0(X,\varphi^*(\widetilde{M})) \xraq{\cong} H^0(X, \cF) \, =: \, M_0
\end{equation*}
is the identity.
\epp 

In this Proposition, we are given a proper morphism $\pi : X \ra Y$ together with an ample invertible sheaf $\cL$. However, in examples of $\bZ$-graded rings coming from flips and flops, instead of an invertible sheaf $\cL$, we are often given a relatively ample Weil divisor $D$ that is only assumed to be $\bQ$-Cartier, \ie, $\cO(dD)$ is an invertible sheaf, but $\cO(D)$ itself may not be. 
We formulate this by considering the following situation ({\it cf}. \cite[Proposition 3.5]{Yeu20a}):
\begin{equation}  \label{sheaf_of_alg_on_X}
\parbox{40em}{
On a scheme $X$ proper over another Noetherian scheme $Y$, there is a quasi-coherent sheaf of $\bZ$-graded algebras $\bigoplus_{i \in \bZ} \cO(i)$ such that each $\cO(i)$ is coherent, and there exists some positive integer $d>0$ such that\\
(1) the sheaf $\cO(d)$ is an invertible sheaf on $X$ ample over $Y$; \\
(2) for each $i,j \in \bZ$, the multiplication map $\cO(di) \otimes_{\cO_{X}} \cO(j) \ra \cO(di+j)$ is an isomorphism.}
\end{equation}

For notational simplicity, we consider here the case when $Y$ is affine, say $Y = \Spec\, R$. 
Let $B$ be the graded algebra over $R$ given by $B_i := H^0(X,\cO(i))$. Then we have the following

\bpp  \label{X_isom_Proj_2}
The $\bN$-graded ring $B_{\geq 0}$ is Noetherian.
Moreover, there is a canonical isomorphism $\varphi: X \xra{\cong} \Proj^+(B)$ over $\Spec\, R$, together with a canonical isomorphism  $\bigoplus_{i \in \bZ} \varphi^*(\widetilde{B(i)}) \xra{\cong} \bigoplus_{i \in \bZ} \cO(i)$ of sheaves of $\bZ$-graded algebras on $X$, such that the composition
\begin{equation*}
B_i \xraq{\eqref{M_to_global_tilde_M}} H^0(\Proj^+(B), \widetilde{B(i)}) \, = \,  H^0(X,\varphi^*(\widetilde{B(i)})) \xraq{\cong} H^0(X, \cO(i)) \, =: \, B_i
\end{equation*}
is the identity. As a result, $B$ is positively $\tfrac{1}{d}$-Cartier.

Furthermore, for any quasi-coherent sheaf $\cF$ on $X$, let $M \in \Gr(B)$ be the graded $B$-module given by $M_i := H^0(X, \cF \otimes \cO(i))$, then there is a canonical isomorphism $\varphi^*( \widetilde{M} ) \xra{\cong} \cF$, such that the composition 
\begin{equation*}
M_0 \xraq{\eqref{M_to_global_tilde_M}} H^0(\Proj^+(B), \widetilde{M}) \, = \,  H^0(X,\varphi^*(\widetilde{M})) \xraq{\cong} H^0(X, \cF) \, =: \, M_0
\end{equation*}
is the identity.
\epp

\bpf
All these statements follow by applying Proposition \ref{X_isom_Proj_1} to the ample invertible sheaf $\cL := \cO(d)$, since we have $\psi : \Proj^+(B^{(d)}) \xra{\cong} \Proj^+(B)$, together with a canonical isomorphism $\psi^*(\widetilde{M}) \cong \widetilde{M^{(d)}}$ for each $M \in \Gr(B)$.
\epf


Now we get back to the situation of a $\bZ$-graded ring $A$ such that $A_{\geq 0}$ is Noetherian, so that it is positively $\tfrac{1}{d}$-Cartier for some $d>0$.
The maps \eqref{M_to_global_tilde_M}, taken for various choices of $M$, assemble to give many other maps.
For example, if we let $\widetilde{M}(i) := \widetilde{M} \otimes_{\cO_{X^+}} \widetilde{A(i)}$, then these maps give rise to canonical graded $A$-module structures on
\begin{equation*}
N' \, := \, \bigoplus_{i \in \bZ} H^0(X^+ , \widetilde{M(i)}) \qquad \text{and} \qquad 
N \, := \, \bigoplus_{i \in \bZ} H^0(X^+ , \widetilde{M}(i))
\end{equation*}
and there are canonical maps $M \raq N' \laq N $ of graded $A$-modules. 
These map satisfies the following properties (see, e.g., \cite[0B5R]{Sta} for the case when $M$ is finitely generated, which implies the general case by taking a directed colimit):

\bpp  \label{MNNp}
Assume that $A_{\geq 0}$ is Noetherian and $A$ is positively $\tfrac{1}{d}$-Cartier, then
\begin{enumerate}
	\item For any $i \in \bZ$, the map $N_{di} \ra N'_{di}$ is an isomorphism.
	\item Both the kernel and cokernel of $M \ra N'$ are $(I^+)^{\infty}$-torsion in the sense of Definition \ref{I_infty_torsion_def}.
	\item If $M_{\geq c}$ is finitely generated over $A_{\geq 0}$ for some $c \in \bZ$, then so are $N_{\geq c'}$ and $N'_{\geq c'}$, for any $c' \in \bZ$.
\end{enumerate}
\epp

\bpf
(1) is obvious. (3) follows from \cite[Tag 0B5R, statements (1),(2)]{Sta}. Thus it suffices to show (2). For this, one can either modify the proof of \cite[Tag 0B5R, statement (5)]{Sta}, or simply notice that $N' = H^0(\Ce_{I^+}(M))$, so that the exact sequence \eqref{Gam_M_Ce_four_terms} gives
\begin{equation*}
\ker(M\ra N') = H^0(\RGam_{I^+}(M)) \qquad \text{and} \qquad 
\coker(M \ra N') = H^1(\RGam_{I^+}(M))
\end{equation*}
both of which are $(I^+)^{\infty}$-torsion (see Theorem \ref{RGam_right_adj}).
\epf


Propositions \ref{X_isom_Proj_2} and \ref{MNNp} can be combined to describe the abelain category of quasi-coherent sheaves on $\Proj^+(A)$ in terms of a Serre quotient (see Theorem \ref{Serre_equiv_thm} below).
We start with the following 

\bdf  \label{uple_torsion_def}
A homogeneous element $x\in M$ in a graded module $M \in \GrA$ is said to be \emph{positively torsion} if it is $(I^+)^{\infty}$-torsion in the sense of Definition \ref{I_infty_torsion_def}. The graded module $M$ is said to be \emph{positively torsion} if every homogeneous element in it is positively torsion. Denote by $\Torp(A) \subset \GrA$ the full subcategory of positively torsion modules. 

For any $m > 0$, denote by $\Torp_{(m)}(A) \subset \Gr(A)$ the full subcategory consisting of graded modules $M \in \GrA$ whose $m$-uple component $M^{(m)} \in \Gr(A^{(m)})$ is in $\Torp(A^{(m)})$. 



Denote by $I^+_{(m)} \subset A^{(m)}$ the graded ideal $I^+_{(m)} := (A^{(m)})_{>0} \cdot A^{(m)}$, and similarly $I^-_{(m)} := (A^{(m)})_{<0} \cdot A^{(m)}$.
\edf

If $A_{\geq 0}$ is Noetherian, then so is $(A^{(m)})_{\geq 0}$. The characterization of positively torsion elements in Lemma \ref{IIT_element_equiv} then shows that $\Torp(A) \subset \Torp_{(m)}(A)$,
so that the Serre subcategory $\Torp_{(m)}(A)$ descends to a Serre subcategory 
$\overline{\Torp_{(m)}(A)} \subset \QpGr(A)$ of the Serre quotient
\begin{equation*}
\QpGr(A) \, := \, \GrA / \Torp(A)
\end{equation*}
and we have a canonical equivalence of abelian categories
\begin{equation*}
\QpmGr(A) \, := \, \GrA / \Torp_{(m)}(A) \, \simeq \, \QpGr(A) / \overline{\Torp_{(m)}(A)}
\end{equation*}

Consider the two functors
\begin{equation}  \label{Serre_equiv_functors_1}
\begin{split}
(-)^{\sim} \,&:\, \GrA \raq \QCoh(X^+)\, , \qquad M \, \mapsto \, \widetilde{M} \\
\!\,^0 \! \cL^{+} \,&:\, \QCoh(X^+) \raq \GrA \, , \qquad \!\,^0 \! \cL^{+}(\cF)_i \, := \, H^0(X^+, \cF \otimes \widetilde{A(i)})
\end{split}
\end{equation}

Suppose that $A_{\geq 0}$ is Noetherian and (postively) $\tfrac{1}{d}$-Cartier, then clearly the functor $(-)^{\sim}$ vanishes on the Serre subcategory $\Torp_{(d)}(A)$, and hence descends to functors 
\begin{equation}  \label{Serre_equiv_functors_2}
\begin{tikzcd}
\!\,^0 \! \cL^{+} \,:\, \QCoh(X^+) \ar[r, shift left]
& \QpdGr(A) \, : \, (-)^{\sim}  \ar[l, shift left]
\end{tikzcd}
\end{equation}
then Propositions \ref{X_isom_Proj_2} and \ref{MNNp} can be combined to show the following result (see also the similar statement in \cite[Tag 0BXF]{Sta} for coherent, instead of quasi-coherent, sheaves on $\Proj^+(A)$, in the case when $A$ is concentrated in non-negative weights).

\bthm [Serre's equivalence]  \label{Serre_equiv_thm}
If $A_{\geq 0}$ is Noetherian and (postively) $\tfrac{1}{d}$-Cartier, then the functors in \eqref{Serre_equiv_functors_2} are inverse equivalences. 
\ethm

\bpf
By Proposition \ref{X_isom_Proj_2}, we see that the composition $(-)^{\sim} \circ \!\,^0 \! \cL^{+} $ is isomorphic to the identity functor on $\QCoh(X^+)$.

The composition $\!\,^0 \! \cL^{+} \circ (-)^{\sim}$ assigns to each $M \in \GrA$ the graded module $N$ in the statement of Poposition \ref{MNNp}. By this Proposition, $M$ and $N$ are linked by a zig-zag of maps $M \ra N' \la N$ in $\GrA$, both of which descend to isomorphisms in $\QpmGr(A)$. 
\epf

Recall that the scheme $X^+ = \Proj^+(A)$ depends only on the ``tail'' of $A$. 
While $\QCoh(X^+)$ is equivalent to $\QpGr(A)$ in the usual way only if $A$ is positively Cartier, it turns out that in general, the latter also depends only on the ``tail'' of $A$.
More precisely, consider a map $f : A \ra B$ of $\bZ$-graded rings, which then induces an adjunction 
\begin{equation} \label{Gr_adj_2}
\begin{tikzcd}
- \otimes_A B \, : \, \Gr(A) \ar[r, shift left] & \Gr(B) \, : \, (-)_A \ar[l, shift left]
\end{tikzcd}
\end{equation}
In general, the right adjoint $(-)_A$ preserves the subcategories of positively torsion modules, and hence induces a functor
$
(-)_A :  \QpGr(B) \ra \QpGr(A)
$.
However, the functor $-\otimes_A B$ may not preserve torsion modules. There is however a special case where it does:

\bpp[\cite{AZ94}, Proposition 2.5]
Suppose that the map $f : A \ra B$ of $\bZ$-graded rings induces isomorphisms $f : A_n \xra{\cong} B_n$ for sufficiently large $n$, then  
the functor $-\otimes_A B : \Gr(A) \ra \Gr(B)$ sends $\Torp(A)$ to $\Torp(B)$, and induces an equivalence $\QpGr(A) \xra{\simeq} \QpGr(B)$.
\epp

In particular, if we apply this to the inclusion map $f : A_{\geq 0} \ra A$, then we have the following
\bcor   \label{QGrA_A_geq0}
The adjunction \eqref{Gr_adj_2} for the map $f : A_{\geq 0} \ra A$ of graded rings decend to give an equivalence $\QpGr(A) \simeq \QpGr(A_{\geq 0})$.
\ecor


\brm  \label{stacky_Proj}
When $A$ is not Cartier, the category $\QpGr(A)$ also admit a description in terms of quasi-coherent sheaves on a stacky projective space. 
Namely, one can show that the map of $\bG_m$-equivariant schemes $\Spec(A) \ra \Spec(A_{\geq 0})$ induces an isomorphism on the $\bG_m$-invariant open subschemes
\begin{equation*}
\Spec(A) \setminus \Spec(A/I^+) \xraq{\cong} \Spec(A_{\geq 0}) \setminus \Spec(A_0)
\end{equation*}
Indeed, this follows from the isomorphism $(A_{\geq 0})_f \xra{\cong} A_f$ for $\deg(f) > 0$.
If we denote this $\bG_m$-equivariant scheme by $W^{ss}(+)$, and let $\, \PProj^+(A)$ be the quotient stack $[W^{ss}(+)/\bG_m]$, then we have an equivalence
$
\QpGr(A)  \simeq  \QCoh(\PProj^+(A)) 
$.
See, e.g., \cite[Proposition 2.3]{AKO08} or 
\cite[Example 2.15]{HL15}.
\erm

Now we relate the derived category $\cD(\QpGr(A))$ to the subcategory $\cD_{\IpTR}(\GrA) \subset \cD(\GrA)$ considered in Section \ref{GM_gr_rings_sec}. See Appendix \ref{Dcat_app} for the notation.

\bpp  \label{DITR_DQpGR}
For each $\spadesuit \in \{\, \, \,   ,+,-,b\}$, the Serre subcategory $\Torp \subset \GrA$ is $\Dsuit$-localizing. The functor ${\bm R}\phi_* : \cD(\QpGr(A)) \ra \cD(\GrA)$ is fully faithful, and there is a semi-orthogonal decomposition 
\begin{equation*}
\Dsuit(\GrA) \, = \, \langle \,  {\bm R}\phi_* ( \Dsuit(\QpGr(A))) \,   , \, \Dsuit_{\Torp}(\GrA) \, \rangle
\end{equation*}
The composition ${\bm R}\phi_* \circ \phi^* : \Dsuit(\GrA) \ra \Dsuit(\GrA) $ is given by $\Ce_{I^+} : \Dsuit(\GrA) \ra \Dsuit(\GrA)$.

As a result, there is an exact equivalence
\begin{equation}  \label{QpGrA_IpTR_equiv}
\begin{tikzcd}
\phi^* \, : \, \Dsuit_{\IpTR}(\GrA) \ar[r, shift left]
&  \Dsuit(\QpGr(A)) \, : \, {\bm R}\phi_* \ar[l, shift left]
\end{tikzcd}
\end{equation}
\epp

\bpf
The Serre subcategory $\Torp \subset \GrA$ is clearly part of a torsion theory on $\GrA$.
By Corollary \ref{Groth_torsion_localizing}, it is therefore localizing. 
By Proposition \ref{Groth_loc_subcat}, it is moreover $\cD$-localizing, so that Proposition \ref{Serre_SOD_1} gives the claimed semi-orthogonal decomposition for the case when  $\spadesuit$ is the empty symbol. 
Thus, the essential image of ${\bm R}\phi_* : \cD(\QpGr(A)) \ra \cD(\GrA)$ is equal to the right orthogonal of $\cD_{\Torp}(\GrA)$, which is therefore equal to $\cD_{\IpTR}(\GrA)$ by the semi-orthogonal decomposition in the first row of \eqref{GM_two_SODs}.
This also shows that ${\bm R}\phi_* \circ \phi^* \cong \Ce_{I^+}$ as endofunctors on $ \cD(\GrA)$.
Since the functor $\Ce_{I^+}$ has finite cohomological dimension, it restricts to each of the subcategories $\Dsuit(\GrA)$. Equivalently, this means that the functor ${\bm R}\phi_* : \cD(\QpGr(A)) \ra \cD(\GrA)$ also restrict to give functors on each of the subcategories $\Dsuit(-)$.
\epf

Now we consider how the equivalence \eqref{QpGrA_IpTR_equiv} restricts to give equivalences on subcategories with coherent cohomologies. The appropriate notions are given in the following
\bdf  \label{IpTR_coh_def}
Denote by  $\gr(A) \subset \GrA$ the full subcategory of finitely generated graded modules.
Let $\qpgr(A) \subset \QpGr(A)$ be the essentially image of $\gr(A)$ under $\phi^* : \GrA \ra \QpGr(A)$.
For each $\spadesuit \in \{\, \, \, ,+,-,b\}$, 
\begin{enumerate}
	\item let $\Dsuit_{\coh}(\QpGr(A)) \subset \Dsuit(\QpGr(A))$ be the full subcategory consisting of complexes whose cohomology lies in $\qpgr(A)$; and
	\item let $\Dsuit_{{\rm coh}(\IpTR)}(\GrA) \subset \Dsuit_{\IpTR}(\GrA)$ be the essential image of $\Dsuit_{{\rm coh}}(\GrA)$ under the functor $\Ce_{I^+} : \Dsuit(\GrA) \ra \Dsuit_{\IpTR}(GrA)$.
\end{enumerate}
\edf

\brm 
\begin{enumerate}
	\item Under the equivalence $\QpGr(A) \simeq \QCoh(\PProj^+(A))$ in Remark \ref{stacky_Proj}, the subcategory $\Dbcoh(\QpGr(A)) \subset \cD(\QpGr(A))$ corresponds
	to the subcategory $\Dbcoh(\PProj^+(A)) \subset \cD(\QCoh(\PProj^+(A)))$, which is often what one is mostly interested in.
	\item The reader is cautioned that in general, $\Dsuit_{{\rm coh}(\IpTR)}(\GrA) \neq \Dsuit_{{\rm coh}}(\GrA) \cap \Dsuit_{\IpTR}(\GrA)$, because the functor $\Ce_{I^+}$ does not preserve $\Dsuit_{\coh}(\GrA)$.
\end{enumerate}
\erm

\blm
For $\spadesuit \in \{-,b\}$, the full subcategory $\Dsuit_{\coh}(\QpGr(A)) \subset \Dsuit(\QpGr(A))$ is equal to the essential image of $\Dsuit_{\coh}(\GrA)$ under the functor $\phi^* : \Dsuit(\GrA) \ra \Dsuit(\QpGr(A))$.
\elm

\brm
We believe that this lemma holds for all cases  $\spadesuit \in \{\, \, \,   ,+,-,b\}$. For example, it seems that one could modify the proof of \cite[Tag 06XL]{Sta} by carefully keeping track of finite generatedness. However, we will be content to give here a more artificial proof by combining several known results.
\erm

\bpf It is clear that the essential image of $\Dsuit_{\coh}(\GrA)$ under $\phi^*$ lies in $\Dsuit_{\coh}(\QpGr(A))$. Conversely, given $\cM \in \Dsuit_{\coh}(\QpGr(A))$, we take ${\bm R}\phi_*(\cM) \in \Dsuit(\GrA)$. In general, this complex will not have finitely generated cohomologies. However, we claim that its restriction to non-negative weights, thought of as an object ${\bm R}\phi_*(\cM)_{\geq 0} \in \Dsuit(\Gr(A_{\geq 0}))$ have cohomologies that are finitely generated over $A_{\geq 0}$. Since the functor ${\bm R}\phi_* \circ \phi^*$ is identified with $\Ce_{I^+}$ in Proposition \ref{DITR_DQpGR}, we see that ${\bm R}\phi_*$ has finite cohomological dimension, so that it suffices to prove this claim for a module $\cM \in \qpgr(A)$, say $\cM = \phi^*(M)$ for $M \in \gr(A)$. Under the identification \eqref{Ce_i_RGam}, this claim then follows from Serre's vanishing, which asserts that $H^j(\Ce_{I^+}(M)_i) = 0$ for all $j >0$ and $i \gg 0$, as well as from Proposition \ref{MNNp}(c), noticing that $N' = H^0(\Ce_{I^+}(M))$. Now, back to the general situation of any ${\bm R}\phi_*(\cM) \in \Dsuit(\GrA)$ for $\spadesuit \in \{-,b\}$. Since ${\bm R}\phi_*(\cM)_{\geq 0} \in \Dmcoh(\Gr(A_{\geq 0}))$, we have $P' := ({\bm R}\phi_*(\cM)_{\geq 0}) \otimes_{A_{\geq 0}}^{{\bm L}} A \in \Dmcoh(\GrA)$ as well. By Corollary \ref{QGrA_A_geq0}, we see that the images of $P'$ and ${\bm R}\phi_*(\cM)$ under $\phi^*$ are isomorphic. This proves the lemma for $\spadesuit = -$. For $\spadesuit = b$, simply take a good truncation $P = \tau^{\geq m}(P')$ for some $m \ll 0$.
\epf


\bcor  \label{coh_IpTR_QpGr}
For $\spadesuit \in \{-,b\}$, the equivalence \eqref{QpGrA_IpTR_equiv} restricts to give an exact equivalence 
\begin{equation*} 
\begin{tikzcd}
\phi^* \, : \, \Dsuit_{\coh(\IpTR)}(\GrA) \ar[r, shift left]
&  \Dsuit_{\coh}(\QpGr(A)) \, : \, {\bm R}\phi_* \ar[l, shift left]
\end{tikzcd}
\end{equation*}
\ecor


\section{Grothendieck duality and Greenlees-May duality}  \label{GGM_sec}

For a quasi-coherent sheaf  $\cF \in \QCoh(X^+)$, the graded module $\!\,^0 \! \cL^{+}(\cF) \in \GrA$ defined in \eqref{Serre_equiv_functors_1} is usually regarded as the associated graded module. There is a different version $\!\,^0 \! \cR^{+}(\cF)$ of the associated graded module, which coincides with $\!\,^0 \! \cL^{+}(\cF)$ in weights $i \in d\bZ$ if $A$ is $\tfrac{1}{d}$-Cartier (see Remark \ref{cR_Serre_equiv_remark} below). Indeed, consider the functor 
\begin{equation}  \label{Serre_R_functor}
\!\,^0 \! \cR^{+} \,:\, \QCoh(X^+) \raq \GrA \, , \qquad \, \!\,^0 \! \cR^{+}(\cF)_i \, := \, \Hom_{\cO_{X^+}}(\widetilde{A(-i)}, \cF)
\end{equation}

This is right adjoint to the associated sheaf functor:
\begin{equation}  \label{Serre_equiv_functors_4}
\begin{tikzcd}
(-)^{\sim} \,:\, \GrA  \ar[r, shift left]
& \QCoh(X^+)  \, : \, \!\,^0 \! \cR^{+}   \ar[l, shift left]
\end{tikzcd}
\end{equation}
Indeed, by the adjoint functor theorem, $(-)^{\sim}$ must have a right adjoint. By investigating the adjunction isomorphism on the objects $A(-i) \in \GrA$, we see that this right adjoint must be given by \eqref{Serre_R_functor}.

For each $M \in \GrA$, there is clearly a canonical isomorphism 
\begin{equation}  \label{H0_0CeI}
\!\,^0\Ce_{I^+}(M)_i \, \cong \, H^0(X^+, \widetilde{M(i)}) 
\end{equation}
where the left hand side is defined as in Corollary \ref{Cech_der_func}. Hence, for any graded module $M \in \GrA$, there is a canonical map 
\begin{equation}  \label{Serre_Ce_R}
 \!\,^0\Ce_{I^+}(M) \raq \!\,^0 \! \cR^{+} (\widetilde{M})
\end{equation}
whose weight $i$ component is the map 
\begin{equation}  \label{Serre_R_functor_adj_map}
H^0(X^+,\widetilde{M(i)}) \raq \Hom_{\cO_{X^+}}(\widetilde{A(-i)}, \widetilde{M}) 
\end{equation}
induced by the multiplcation map $\widetilde{M(i)} \otimes_{\cO_{X^+}} \widetilde{A(-i)} \ra \widetilde{M}$. 
In particular, the map \eqref{Serre_Ce_R} is an isomorphism in weights $i \in d\bZ$ if $A$ is $\tfrac{1}{d}$-Cartier.

\brm  \label{cR_Serre_equiv_remark}
One usually takes the functor $\!\,^0 \! \cL^{+}$ to define Serre's equivalence, as in \eqref{Serre_equiv_functors_2}. An alternative choice is the functor $\!\,^0 \! \cR^{+}$. Indeed, these two functors are naturally isomorphic after passing to $\QpdGr(A)$.
\erm

We now consider the derived versions of \eqref{Serre_R_functor}, \eqref{Serre_equiv_functors_4}, \eqref{H0_0CeI}, \eqref{Serre_Ce_R}, and \eqref{Serre_R_functor_adj_map}. First, since $\QCoh(X^+)$ is a Grothendieck category, it has enough K-injectives, so that \eqref{Serre_R_functor} has a right derived functor
$\cR^{+} : \cD(\QCoh(X^+)) \ra \cD(\GrA)$, whose weight $i$ component is given by
\begin{equation} \label{RFi_RHom}
\cR^{+}(\cF)_i \, \cong \, \RHom_{\cO_{X^+}}( \widetilde{A(-i)} , \cF )
\end{equation}

There is still an adjunction at the level of derievd categories
\begin{equation} \label{m_uple_adj_derived_1}
\begin{tikzcd}
(-)^{\sim} \, : \, \cD(\GrA) \ar[r, shift left] 
& \cD(\QCoh(X^+)) \, : \,  \cR^{+}   \ar[l, shift left] 
\end{tikzcd}
\end{equation}

To obtain the analogue of \eqref{H0_0CeI}, we give the following graded analogue of \cite[Tag 09T2]{Sta}:
\blm  \label{inj_flasque}
If $A$ is Noetherian, then for any injective object $M \in \GrA$, the associated sheaf $\widetilde{M} \in \QCoh(X^+)$ is flasque.
\elm

\bpf
For any closed subset $Z \subset \Proj^+(A)$, say defined by a graded ideal $I\subset I^+$, let $U = \Proj^+(A) \setminus Z$, then, as in \cite[Tag 01YB]{Sta}, any section $s\in \Gamma(U,\widetilde{M})$ is represented by an element in $\Hom_{\GrA}(I^n,M)$, for some $n \geq 0$. Since $M$ is injective, it extends to an element in $\Hom_{\GrA}(A,M) = M_0$. 
\epf

Thus, combining Corollary \ref{Cech_der_func} and Lemma \ref{inj_flasque}, we see that for each $i \in \bZ$, there is a canonical isomorphism in $\cD(A_0)$:
\begin{equation}  \label{Ce_i_RGam}
\Ce_{I^+}(M)_i \, \cong \, {\bm R}\pi^+_*( \widetilde{M(i)} )
\end{equation}
where ${\bm R}\pi^+_*$ is the derived pushforward functor for $\pi^+ : \Proj^+(A) \ra \Spec(A_0)$.
This isomorphism is, of course, the usual statement that sheaf cohomology can be computed by the \v{C}ech complex.

By replacing \eqref{Serre_Ce_R} by its derived functors%
\footnote{More precisely, by Corollary \ref{Cech_der_func} again, the derived version of \eqref{Serre_Ce_R} gives a natural transformation $\Ce_{I^+} \ra {\bm R}(\!\,^0 \! \cR^{+} \circ (-)^{\sim})$, which is then post-composed with the canonical natural transformation ${\bm R}(F \circ G) \ra {\bm R} F \circ {\bm R} G$ for $F = \!\,^0 \! \cR^{+}$ and $G = (-)^{\sim}$. }, 
we have a canonical map in $\cD(\GrA)$
\begin{equation}  \label{Cech_R_functor}
\Ce_{I^+}(M) \raq \cR^+(\widetilde{M})
\end{equation}
so that under the identification \eqref{Ce_i_RGam}, the weight $i$ component of \eqref{Cech_R_functor} is given by the canonical map 
\begin{equation*}
{\bm R}\pi^+_*( \widetilde{M(i)} ) \raq \RHom_{X^+}(\widetilde{A(-i)}, \widetilde{M})
\end{equation*}

The identification \eqref{Ce_i_RGam} also allows us to show the following two results:
\blm  \label{RGam_Dbcoh}
If $A$ is Noetherian, then for any $M \in \Dbcoh(\GrA)$, we have $\Ce_{I^+}(M)_i \in \Dbcoh(A_0)$ and $\RGam_{I^+}(M)_i \in \Dbcoh(A_0)$ for each weight $i \in \bZ$.
\elm

\bpf
In view of \eqref{Ce_i_RGam}, this follows from the statement that the derived pushforward functor ${\bm R}\pi^+_*$ for the projective morphism $\pi^+$ preserves $\Dbcoh(-)$.
\epf

\blm  \label{local_cohom_weight_bounded}
If $A$ is Noetherian, then for any $M \in \Dbcoh(\GrA)$, there exists $c^+, c^- \in \bZ$ such that 
$\RGam_{I^+}(M)_i \simeq 0$ for all $i > c^+$ and $\RGam_{I^-}(M)_i \simeq 0$ for all $i < c^-$.
\elm

\bpf
Clearly, it suffices to prove the lemma for a finitely generated graded module $M$.
By (2) and (3) of Proposition \ref{MNNp}, we see that $M_i \ra N'_i = H^0(\Ce_{I^+}(M)_i)$ is an isomorphism for $i \geq c_0^+$ for some $c_0^+ \in \bZ$. Thus, for $i \geq c_0^+$, we have $\RGam_{I^+}(M)_i \simeq 0$ if and only if ${\bm R}^j \pi^+_*( \widetilde{M(i)}) = 0$ for all $j > 0$.
Since $A$ is Noetherian, it is $\tfrac{1}{d}$-Cartier for some $d>0$. Then the sequence $\widetilde{M(i)}, \widetilde{M(i+d)}, \widetilde{M(i+2d)},\ldots$ of coherent sheaves on $X^+$ is a sequence of twist by an ample invertible sheaf, and hence must eventually have zero higher cohomology. Apply this for $i = 0,\ldots,d-1$ in order to find $c^+$.
The integer $c^-$ can also be found in a similar way, by considering the sequence 
$\widetilde{A(i)}, \widetilde{A(i-d)}, \widetilde{A(i-2d)},\ldots$ on $X^-$.
\epf

From now on, we assume that $A$ is Noetherian without any more explicit mention.
First, we will relate Greenlees-May duality on the graded ring $A$ and Grothendieck duality for the projective morphism $\pi^+ : X^+ \ra \Spec(R)$, where we have written $R := A_0$. 
%
%
Consider the following composition of adjunctions
\begin{equation}  \label{GGM_comp_adj_1}
\begin{tikzcd}
\cD(\GrA)  \ar[r, shift left, "\Ce_{I^+}"]
& \cD_{\IpTR}(\GrA) \ar[l, shift left, "\iota"]  \ar[r, shift left, "\iota"]
& \cD(\GrA)  \ar[l, shift left, "\cE_{I^+}"]  \ar[rr, shift left, "(-)_0"]
& & \cD(\Mod(R)) \ar[l l, shift left, "\RHomcom_R(A\text{,}-)"]
\end{tikzcd}
\end{equation}

A crucial observation is that, by \eqref{Ce_i_RGam}, the composition of all the right-pointing arrows is naturally isomorphic to the functor $M \mapsto {\bm R}\pi^+_*(\widetilde{M})$.
This allows us to compare \eqref{GGM_comp_adj_1} with the following composition of adjunctions:
\begin{equation}  \label{GGM_comp_adj_2}
\begin{tikzcd}
\cD(\GrA)  \ar[r, shift left, "(-)^{\sim}"]
& \cD(\QCoh(X^+)) \ar[l, shift left, "\cR^+"]  \ar[r, shift left, "{\bm R}\pi^+_*"]
& \cD(\Mod(R)) \ar[l, shift left, "(\pi^+)^!"]
\end{tikzcd}
\end{equation}
where the adjunction pair on the left is given by \eqref{m_uple_adj_derived_1}. This gives the following

\bthm  \label{GGM_thm}
For each $L \in \cD(\Mod(R))$, there is a canonical isomorphism in $\cD(\GrA)$:
\begin{equation}  \label{cR_pi_shriek_isom_1}
\cR^+ ((\pi^+)^! ( L ))  \, \cong \,  \cE_{I^+}( \RHomcom_R(A,L) ) 
\, \cong \, \RHomcom_R(\Ce_{I^+}(A),L)
\end{equation}
Moreover, under the isomorphisms \eqref{RFi_RHom} and \eqref{Ce_i_RGam}, the weight $i$ component of this isomorphism can be described by the commutativity of the following diagram in $\cD(R)$:
\begin{equation}  \label{cR_pi_shriek_isom_2}
\begin{tikzcd}
\cR^+ ((\pi^+)^! ( L ))_i \ar[r, "\eqref{cR_pi_shriek_isom_1}", "\cong"'] \ar[d, "\eqref{RFi_RHom}"', "\cong"]  &   \RHomcom_R(\Ce_{I^+}(A)_{-i},L) \ar[d, "\eqref{Ce_i_RGam}", "\cong"'] \\
\RHom_{\cO_{X^+}}( \widetilde{A(-i)} , (\pi^+)^!(L) ) \ar[r, "\cong"] & \RHom_R( {\bm R}\pi^+_*(\widetilde{A(-i)}),L)
\end{tikzcd}
\end{equation}
where the horizontal map in the bottom is the local adjunction isomorphism.
\ethm

\bpf
The first isomorphism of \eqref{cR_pi_shriek_isom_1} is obtained by comparing \eqref{GGM_comp_adj_1} and \eqref{GGM_comp_adj_2}. The second isomorphism of \eqref{cR_pi_shriek_isom_1} is a consequence of the general fact $\cE_{I}(M) \cong \RHomcom_A(\CI(A), M)$.

To show that the diagram \eqref{cR_pi_shriek_isom_2} commutes, recall that the vertical map \eqref{RFi_RHom} on the left is, by definition the adjunction isomorphism for the first pair $(-)^{\sim} \dashv \cR^+$ of \eqref{GGM_comp_adj_2} applied to $A(-i)$ and $(\pi^+)^!(L)$; while the bottom horizontal arrow of \eqref{cR_pi_shriek_isom_2} is the adjunction isomorphism for the second pair ${\bm R}\pi^+_* \dashv (\pi^+)^!$ of \eqref{GGM_comp_adj_2}. 
Thus, the map $\cR^+ ((\pi^+)^! ( L ))_i \ra \RHom_R( {\bm R}\pi^+_*(\widetilde{A(-i)}),L)$ obtained via the lower route of \eqref{cR_pi_shriek_isom_2} is precisely the adjunction isomorphism for \eqref{GGM_comp_adj_2}. 
Since the isomorphism \eqref{cR_pi_shriek_isom_1} is defined by identifying both as right adjoint of two functors identified via \eqref{Ce_i_RGam}, the corresponding adjunction isomorphisms also get identified under \eqref{Ce_i_RGam}. This translates precisely to the commutativity of \eqref{cR_pi_shriek_isom_2}.
\epf

\brm  \label{GGM_remark}
One may view the isomorphism \eqref{Ce_i_RGam} as an interpretation of the functor ${\bm R}\pi^+_*$ in terms of the functor $\Ce_{I^+}$. Similarly, one may also regard Theorem \ref{GGM_thm} as an interpretation of the Grothendieck duality functor $(\pi^+)^!$ in terms of the functor $\cE_{I^+}$.
Namely, if we think of $\cR^+$ as the functor that gives the associated graded module (see Remark \ref{cR_Serre_equiv_remark}), then Theorem \ref{GGM_thm} says that the associated graded module of $(\pi^+)^!$ is given by $\RHomcom_R(\Ce_{I^+}(A),-)$.
\erm

\brm
The existence of the map \eqref{Cech_R_functor} and the isomorphism \eqref{cR_pi_shriek_isom_1} in $\cD(\GrA)$ should be intuitively quite clear, as there are canonical maps (resp. isomorphisms) relating each of their weight components. However, to actually prove that these assemble to maps in $\cD(\GrA)$, we have found it necessary to develop the relevant functors from scratch, as we do here.
\erm

\section{The case of non-affine base}  \label{non_affine_sec}




In this section, we provide the formal arguments to extend our previous discussion to the case of non-affine base. More precisely, we work in the following setting:
\begin{equation}  \label{sheaf_A_setting}
\parbox{40em}{$Y$ is a Noetherian separated scheme, and $\cA$ is a quasi-coherent sheaf of Noetherian $\bZ$-graded rings on $Y$, such that $\cA_0$ (and hence every $\cA_i$) is coherent over $\cO_X$.}
\end{equation}

Denote by $\Gr(\cA)$ the category of quasi-coherent graded $\cA$-modules. Then $\Gr(\cA)$ is equivalent to the category $\QCoh_{\bG_m}(W)$ of $\bG_m$-equivariant quasi-coherent sheaves on the relative spectrum $W := \Spec_Y \cA$.
We denote this equivalence by
\begin{equation} \label{sharp_flat_GrA}
\begin{tikzcd}
(-)^{\sharp} \, : \, \Gr(\cA) \ar[r, shift left] 
& \QCoh_{\bG_m}(W) \, : \, (-)^{\flat} \ar[l, shift left]
\end{tikzcd}
\end{equation}
We fix our convention so that
\begin{equation}  \label{shift_twist_conv}
\parbox{40em}{The shift functor $\cM \mapsto \cM(1)$ on $\Gr(\cA)$ corresponds under \eqref{sharp_flat_GrA} to the twist of equivariance structure by the identity character of $\bG_m$}
\end{equation}

For any $\cM, \cN \in \Gr(\cA)$, one can define $\cM \otimes_{\cA} \cN \in \Gr(\cA)$, which coincides with the graded tensor product over each affine open subset $U \subset Y$. Moreover, we have $(\cM \otimes_{\cA} \cN)^{\sharp} \cong \cM^{\sharp} \otimes_{\cO_W} \cN^{\sharp}$, where the right hand side inherits an equivariance structure in the usual way.

Now we define local cohomology complex. Given any quasi-coherent sheaf of graded ideals $\scI \subset \cA$, we say that a quasi-coherent sheaf $\cM \in \Gr(\cA)$ of graded module is \emph{$\scI^{\infty}$-torsion} if $\cM^{\sharp}$ restricts to zero on the complement of $V(\scI) := \Spec_Y(\cA/\scI)$ in $W$. Notice that when $Y$ is affine, this coincides with Definition \ref{I_infty_torsion_def}. 
Denote by $U(\scI) := W \setminus V(\scI)$, and $j : U(\scI) \rinto W$ the inclusion. 
Then the subcategory $\cD_{\cIIT}(\Gr(\cA)) \subset \cD(\Gr(\cA))$ consisting of complexes with $\scI^{\infty}$-torsion cohomology sheaves corresponds to the kernel of $j^* : \cD(\QCoh_{\bG_m}(W)) \ra \cD(\QCoh_{\bG_m}(U(\scI)))$ under the equivalence \eqref{sharp_flat_GrA}. Since the functor $j^*$ has a right adjoint ${\bm R}j_*$ satisfying $j^* \circ {\bm R}j_* \cong \id$, it follows formally that there is a semi-orthogonal decomposition 
\begin{equation}  \label{local_cohom_SOD}
\cD(\Gr(\cA)) \, = \, \langle \, \cD_{\cITR}(\Gr(\cA)) \, , \,  \cD_{\cIIT}(\Gr(\cA)) \, \rangle
\end{equation}
where $ \cD_{\cITR}(\Gr(\cA))$ is the essential image of the fully faithful composition 
\begin{equation*}
\cD_{\cITR}(\Gr(\cA)) \, := \, {\rm Ess\, Im} \, \bigl[ \,
\cD(\QCoh_{\bG_m}(U(\scI))) \xra{ {\bm R}j_* } \cD(\QCoh_{\bG_m}(W)) \xra[\simeq]{(-)^{\flat}} \cD(\Gr(\cA))
\, \bigr]
\end{equation*}

More precisely, for any $\cM \in \cD(\Gr(\cA))$, define $\Ce_{\scI}(\cM) := ({\bm R}j_*(j^*(\cM^{\sharp})))^{\flat} \in \cD(\Gr(\cA))$, which comes with a natural adjunction unit $\eta_{\cM} : \cM \ra \Ce_{\scI}(\cM)$. 
Then $\RGam_{\scI}(\cM) \in \cD_{\cIIT}(\Gr(\cA))$ can be defined by the exact triangle
\begin{equation}  \label{RGam_Ce_tri}
\ldots \raq \RGam_{\scI}(\cM) \xraq{\epsilon_{\cM}} \cM \xraq{\eta_{\cM}} \Ce_{\scI}(\cM) \xraq{\delta_{\cM}} \RGam_{\scI}(\cM)[1] \raq \ldots
\end{equation}
which is the decomposition triangle associated to the semi-orthogonal decomposition \eqref{local_cohom_SOD}.
Over affine open subsets of $Y$, this coincides with the corresponding notions in Section \ref{GM_gr_rings_sec}:
\blm  \label{RGam_affine_open}
The exact triangle \eqref{RGam_Ce_tri} restricts to \eqref{RGam_Ce_seq}
over each affine open subset $U \subset Y$.
\elm

\bpf
Clearly, the association of \eqref{RGam_Ce_tri} to $\cM \in \Gr(\cA)$ is local on the base $Y$, so that it suffices to assume that $Y=U$ is affine. In this case, 
the notion of $\scI^{\infty}$-torsion modules clearly coincides with Definition \ref{I_infty_torsion_def}. 
Since $\RGam_I : \cD(\GrA) \ra \cD_{\IIT}(\GrA)$ was known to be right adjoint to the inclusion by Theorem \ref{RGam_right_adj}, it must coincide with $\RGam_{\scI}$.
\epf

One can use it to show that the isomorphisms $\RGam_I(M) \cong \RGam_I(A) \otimes_A^{{\bm L}} M$ and $\Ce_I(M) \cong \Ce_I(A) \otimes_A^{{\bm L}} M$  in $\cD(\GrA)$ still hold in the present non-affine case:
\bcor  \label{RGam_Ce_tensor}
For any $\cM \in \cD(\Gr(\cA))$, there are canonical isomorphisms $\RGam_{\scI}(\cM) \cong \RGam_{\scI}(\cA) \otimes_{\cA}^{{\bm L}} \cM$ and $\Ce_{\scI}(\cM) \cong \Ce_{\scI}(\cA) \otimes_{\cA}^{{\bm L}} \cM$ in $\cD(\Gr(\cA))$. 
\ecor

%

We now turn to local homology complex in the non-affine setting. Already in the ungraded case, this notion has a technical subtlety arising from the fact that the internal Hom between quasi-coherent sheaves may not be quasi-coherent (see, e.g., \cite[Remark (0.4)]{ATJLJ97} for a discussion). The notion of (quasi-)coherator from \cite[Appendix B]{TT90} is therefore highly relevant in this discussion. We recall this notion now.

On a quasi-compact separated scheme $X$, the inclusion $\iota : \QCoh(X) \ra \Mod(\cO_X)$ has a right adjoint, called the \emph{(quasi-)coherator} $Q_X : \Mod(\cO_X) \ra \QCoh(X)$. The derived functor ${\bm R}Q_X : \cD(X) \ra \cD(\QCoh(X))$ then restricts to an equivalence ${\bm R}Q_X : \Dqcoh(X) \ra \cD(\QCoh(X))$, which is inverse to the inclusion
(see, e.g., \cite[Tag 08DB]{Sta}). 
In this paper, we work with $\cD(\QCoh(X))$ instead of $\Dqcoh(X)$. More precisely, we work with the following two conventions:
\begin{conv}  \label{Dcat_conv1}
	On a quasi-compact separated scheme $X$, let $\cD(X) = \cD(\Mod(\cO_X))$ be the derived category of all sheaves of $\cO_X$-modules.
	
	For any $\cF,\cG \in \Mod(\cO_X)$, denote by $\cHom^{\clubsuit}_{\cO_X}(\cF,\cG) \in \Mod(\cO_X)$  the internal Hom object in $\Mod(\cO_X)$, defined by $\cHom^{\clubsuit}_{\cO_X}(\cF,\cG)(U) = \Hom_{\cO_U}(\cF|_U, \cG|_U)$.
	
	For any $\cF,\cG \in \QCoh(X)$, denote by $\cHom_{\cO_X}(\cF,\cG) \in \QCoh(X)$ the internal Hom object in $\QCoh(X)$, defined by $\cHom_{\cO_X}(\cF,\cG) = Q_X( \cHom^{\clubsuit}_{\cO_X}(\cF,\cG) )$.
	
	Both of these internal Hom functors have derived functors
	\begin{equation*}
	\begin{split}
	\RcHom^{\clubsuit}_{\cO_X}(-,-) \, &: \, \cD(X)^{\op} \times \cD(X) \raq \cD(X) \\
	\RcHom_{\cO_X}(-,-) \, &: \, \cD(\QCoh(X))^{\op} \times \cD(\QCoh(X)) \raq \cD(\QCoh(X))
	\end{split}
	\end{equation*}
	so that for any $\cF,\cG \in \cD(\QCoh(X))$, we have 
	\begin{equation*}
	\RcHom_{\cO_X}(\cF,\cG) \, \cong \, {\bm R}Q_X( \RcHom^{\clubsuit}_{\cO_X}(\cF,\cG) )
	\end{equation*}
	
	In particular, if $X$ is Noetherian, $\cF \in \Dmcoh(\QCoh(X))$ and $\cG \in \cD^+(\QCoh(X))$, then the canonical map 
	$ \RcHom_{\cO_X}(\cF,\cG) \ra  \RcHom^{\clubsuit}_{\cO_X}(\cF,\cG)$ is an isomorphism in $\cD(X)$ (see, e.g., \cite[Tag 0A6H]{Sta}). As a result, these two will often be implicitly identified.
\end{conv}

\begin{conv}  \label{Dcat_conv2}
	Let $f : X \ra Y$ be a quasi-compact morphism between quasi-compact separated schemes. The functors $({\bm L}f^*, {\bm R}f_*, f^!)$ are regarded as functors between $\cD(\QCoh(X))$ and $\cD(\QCoh(Y))$. Namely, ${\bm L}f^*$ and ${\bm R}f_*$ are the derived functors of the functors $f^*$ and $f_*$ between $\QCoh(X)$ and $\QCoh(Y)$.
	Under the equivalences ${\bm R}Q : \Dqcoh(-) \xra{\simeq} \cD(\QCoh(-))$, the functor ${\bm R}f_*$ coincides with the usual one between $\Dqcoh(X)$ and $\Dqcoh(Y)$ (see, e.g., \cite[Tag 0CRX]{Sta}). By adjunction, the same holds for ${\bm L}f^*$, as well as $f^!$, whenever it is well-defined%
	\footnote{In this paper, we will only encounter the functor $f^!$ for proper $f$. In this case, it is defined as the right adjoint of ${\bm R}f_*$.}.
\end{conv}

To develop local homology in the graded context, a careful discussion of sheafified graded Hom complexes and the graded (quasi-)coherator is in order. 
For the sake of a formal argument, we temporarily introduce the category $\GrMod(\cA)$ of (not necessarily quasi-coherent) sheaves of graded $\cA$-modules. For any $\cM , \cN \in \GrMod(\cA)$, define $\cHomcom^{\clubsuit}_{\cA}(\cM,\cN) \in \GrMod(\cA)$ by 
\begin{equation*}
\cHomcom^{\clubsuit}_{\cA}(\cM,\cN)(U)_i \, := \, \Hom_{\GrMod(\cA|_U)}(\cM|_U , \cN|_U(i))
\end{equation*}

This gives a bifunctor
\begin{equation*}
\cHomcom^{\clubsuit}_{\cA}(-,-) \, : \, \GrMod(\cA)^{\op} \times \GrMod(\cA) \raq \GrMod(\cA)
\end{equation*}

Since we are working over the category $\GrMod(\cA)$ of not necessarily quasi-coherent sheaves, the restriction functor $j^* : \GrMod(\cA) \ra \GrMod(\cA|_U)$, for any open subscheme $j : U \rinto Y$, has an exact left adjoint $j_!$ (see, e.g., \cite[Tag 0797]{Sta}). Thus, the restriction of any K-injective complex is K-injective. As a result, one can use K-injective complexes to define the sheafified derived Hom complex:
\begin{equation}  \label{RHom_clubsuit}
\RcHomcom^{\clubsuit}_{\cA}(-,-) \, : \, \cD(\GrMod(\cA))^{\op} \times \cD(\GrMod(\cA)) \raq \cD(\GrMod(\cA))
\end{equation}

This is the internal Hom object in $\cD(\GrMod(\cA))$, with respect to the monoidal product $-\otimes_{\cA}^{{\bm L}} -$. In other words, for any  $\cM,\cN,\cK \in \cD(\GrMod(\cA))$, there is a canonical isomorphism
\begin{equation}  \label{RHom_adj_1}
\Hom_{\cD(\GrMod(\cA))}(\cM \otimes_{\cA}^{{\bm L}} \cN ,\cK) \, \cong \, \Hom_{\cD(\GrMod(\cA))}(\cM, \RcHomcom^{\clubsuit}_{\cA}(\cN,\cK))
\end{equation}

By the adjoint functor theorem,
the inclusion functor $\iota : \Gr(\cA) \ra \GrMod(\cA)$ has a right adjoint $Q_{\cA} : \GrMod(\cA) \ra \Gr(\cA)$, known as the \emph{graded (quasi-)coherator}.
Moreover, this functor can be derived to obtain ${\bm R}Q_{\cA} : \cD(\GrMod(\cA)) \ra \cD(\Gr(\cA))$, which is still right adjoint to the $\iota : \cD(\Gr(\cA)) \ra \cD(\GrMod(\cA))$. 
This allows us to define the bifunctor
\begin{equation}  \label{RGam_DGrA}
\RcHomcom_{\cA}(-,-) \, : \, \cD(\Gr(\cA))^{\op} \times \cD(\Gr(\cA)) \raq \cD(\Gr(\cA))
\end{equation}
by post-composing \eqref{RHom_clubsuit} with the derived (quasi-)coherator: 
\begin{equation*}
\RcHomcom_{\cA}(\cM,\cN) \, := \,  {\bm R}Q_{\cA} \, \RcHomcom^{\clubsuit}_{\cA}(\cM,\cN)
\end{equation*}

Since ${\bm R}Q_{\cA} : \cD(\GrMod(\cA)) \ra \cD(\Gr(\cA))$ is right adjoint to the inclusion, it follows from \eqref{RHom_adj_1} that, for any $\cM,\cN,\cK \in \cD(\Gr(\cA))$, there is a canonical isomorphism
\begin{equation}  \label{RHom_adj_2}
\Hom_{\cD(\Gr(\cA))}(\cM \otimes_{\cA}^{{\bm L}} \cN ,\cK) \, \cong \, \Hom_{\cD(\Gr(\cA))}(\cM, \RcHomcom_{\cA}(\cN,\cK))
\end{equation}
so that $\RcHomcom_{\cA}(-,-)$ is the internal Hom bifunctor in $\cD(\Gr(\cA))$.

In particular, if we define the functors
\begin{equation}  \label{LLam_Ce_functors}
\begin{split}
\LLam_{\scI} \, &: \, \cD(\Gr(\cA)) \ra \cD(\Gr(\cA)) \, \qquad \LLam_{\scI}(\cM) \, := \, \RcHomcom_{\cA}(\RGam_{\scI}(\cA),\cM) \\
\cE_{\scI} \, &: \, \cD(\Gr(\cA)) \ra \cD(\Gr(\cA)) \, \qquad \cE_{\scI}(\cM) \, := \, \RcHomcom_{\cA}(\Ce_{\scI}(\cA),\cM) 
\end{split}
\end{equation}
then by \eqref{RHom_adj_2} and Corollary \ref{RGam_Ce_tensor}, there are canonical adjunctions
\begin{equation}  \label{RGam_LLam_Ce_cE_adj}
\begin{tikzcd}[row sep = 0]
\RGam_{\scI} \, : \,  \cD(\Gr(\cA)) \ar[r,shift left] & \cD(\Gr(\cA)) \, :\,  \LLam_{\scI} \ar[l,shift left] \\
\Ce_{\scI} \, : \, \cD(\Gr(\cA)) \ar[r,shift left] & \cD(\Gr(\cA)) \, : \, \cE_{\scI} \ar[l,shift left]
\end{tikzcd}
\end{equation}

Now we prove the following

\bpp  \label{RQ_A_inverse}
The canonical functor $\iota : \cD(\Gr(\cA)) \ra \cD_{\Gr(\cA)}({\rm GrMod}(\cA))$ is an equivalence, with inverse given by ${\bm R}Q_{\cA}$.
\epp

\bpf
The corresponding statement for $\iota : \cD(\QCoh(Y)) \ra \cD_{\QCoh}(\Mod(\cO_Y))$ is well-known (see, e.g., \cite[Appendix B]{TT90} or \cite[Tag 08DB]{Sta}). One can either adapt this proof to our present case, or formally deduce our statement from that, as follows:

For each $i \in \bZ$, consider the diagram of functors
\begin{equation*}
\begin{tikzcd}
\cD(\Gr(\cA)) \ar[r, "\iota"]
& \cD({\rm GrMod}(\cA)) 
& &  \cD(\Gr(\cA)) \ar[d, "(-)_i"']
& \cD({\rm GrMod}(\cA))  \ar[l, "{\bm R}Q_{\cA} "'] \ar[d, "(-)_i"]\\
\cD(\QCoh(Y)) \ar[r, "\iota"] \ar[u, "-\otimes_{\cO_Y}^{{\bm L}} \cA(-i)"]
& \cD(\Mod(\cO_Y)) \ar[u, "-\otimes_{\cO_Y}^{{\bm L}} \cA(-i)"']
& & \cD(\QCoh(Y)) 
& \cD(\Mod(\cO_Y)) \ar[l, "{\bm R}Q_Y"']
\end{tikzcd}
\end{equation*}
Since the left diagram commutes up to isomorphism of functors, if we take the right adjoints of all the functors involved, we see that the right diagram also commutes up to isomorphism of functors. 
Thus, if we forget about the $\cA$-module structure, then the derived (quasi-)coherator ${\bm R}Q_{\cA}(\cM)$ is simply obtained by applying the derived (quasi-)coherator ${\bm R}Q_Y$ to each weight component $\cM_i$. 
Therefore, the fact that the adjunction unit $\id \Rightarrow {\bm R}Q_{\cA} \circ \iota$ is an isomorphism on $\cD(\Gr(\cA))$; and the adjunction counit $\iota \circ {\bm R}Q_{\cA} \Rightarrow \id$ is an isomorphism on $\cD_{\Gr(\cA)}(\GrMod(\cA))$, follows from the corresponding statements for the adjunction $\iota \dashv {\bm R}Q_Y$.
\epf

\brm
A disadvantage of \eqref{RGam_DGrA}, and hence of the functors \eqref{LLam_Ce_functors}, is that it often fails to be local. Namely, if $U \subset Y$ is an affine open subscheme, and if we let $A := \cA(U)$, then it may happen that $\RcHomcom_{\cA}(\cM,\cN)(U) \not\simeq \RcHomcom_{A}(\cM(U),\cN(U)$.
However, one can show that if $\cM \in \Dmcoh(\Gr(\cA))$ and $\cN \in \cD^+(\Gr(\cA))$, then equality holds. 
\erm

For any $\cM \in \GrMod(\cA)$ and $\cG \in \Mod(\cO_Y)$, define the graded sheafified Hom $\underline{\cHom}^{\clubsuit}_{\cO_Y}( \cM , \cG ) \in \GrMod(\cA)$ by 
\begin{equation*}
\underline{\cHom}^{\clubsuit}_{\cO_Y}( \cM , \cG )(U)_i \, := \, 
\Hom_{\cO_Y}( (\cM|_U)_{-i} , \cG|_U )
\end{equation*}
By taking K-injective representatives of any $\cG \in \cD(\Mod(\cO_Y))$, one can define its derived version%
\footnote{We are implciitly appealing to the argument in the paragraph preceding \eqref{RHom_clubsuit}.}.
This gives a bifunctor
\begin{equation}  \label{RHom_club_A_Y}
\RcHomcom_{\cO_Y}^{\clubsuit}(-,-) \, : \, \cD({\rm GrMod}(\cA))^{\op} \, \times \, \cD(\Mod(\cO_Y)) \raq \cD({\rm GrMod}(\cA))
\end{equation} 
which then allows us to define 
\begin{equation}  \label{RHom_A_Y}
\RcHomcom_{\cO_Y}(-,-) \, : \, \cD(\Gr(\cA))^{\op} \times \cD(\QCoh(Y)) \raq \cD(\Gr(\cA))
\end{equation}
by post-composing \eqref{RHom_club_A_Y} with the derived (quasi-)coherator: 
\begin{equation*}
\RcHomcom_{\cO_Y}(\cM,\cG) \, := \,  {\bm R}Q_{\cA} \, \RcHomcom^{\clubsuit}_{\cO_Y}(\cM,\cG)
\end{equation*}

By the argument in the proof of Proposition \ref{RQ_A_inverse}, we see that the derived (quasi-)coherator commutes with taking weight components. Thus, we have
\begin{equation}  \label{RHom_weight_comp}
\RcHomcom_{\cO_Y}(\cM,\cG)_i \, \simeq \, \RcHom_{\cO_Y}(\cM_{-i},\cG)
\end{equation}

For each $\cM , \cN \in \cD(\Gr(\cA))$ and $\cG \in \cD(\QCoh(Y))$, there is a canonical isomorphism in $\cD(\Gr(\cA))$:
\begin{equation}  \label{RHom_A_Y_adj}
\RcHomcom_{\cO_Y}(\cM \otimes_{\cA}^{{\bm L}} \cN, \cG)
\, \cong \, 
\RcHomcom_{\cA}(\cM, \RcHomcom_{\cO_Y}(\cN,\cG))
\end{equation}


In certain cases of interest, the bifunctor \eqref{RHom_A_Y} is local on the base:
\blm  \label{RHom_Y_qcoh}
Suppose $\cM \in \cD(\Gr(\cA))$ is such that $\cM_i \in \Dmcoh(\QCoh(Y))$ for each weight component $i \in \bZ$, then for any $\cG \in \cD^+(\QCoh(Y))$, the complexes
$\RcHomcom_{\cO_Y}(\cM,\cG)$ and $\RcHomcom_{\cO_Y}^{\clubsuit}(\cM,\cG)$ represent the same object in $\cD({\rm GrMod}(\cA))$. In particular, for any affine open subscheme $U = \Spec \, R \subset Y$, we have
\begin{equation*} 
\RcHomcom_{\cO_Y}(\cM,\cG)(U) \cong \RHomcom_{R}(\cM(U),\cG(U))
\end{equation*}
as objects in $\cD(\Gr(A))$, where $A = \cA(U)$.
\elm


An important special case is when $\cG$ is a dualizing complex $\omega^{\bullet}_Y \in \Dbcoh(\QCoh(Y))$. In this case, $\RcHomcom_{\cO_Y}(\cM,\cG)$ is local on the base when $\cM \in \cD_{\lc}(\Gr(\cA))$, in the sense of the following 
\bdf
A quasi-coherent sheaf $\cM \in \Gr(\cA)$ of graded $\cA$-modules is said to be \emph{locally coherent} if each $\cM_i$ is coherent as a sheaf over $Y$. 
For each $\spadesuit \in \{ \, \, , +,-,b\}$, denote by $\Dsuit_{\lc}(\Gr(\cA))$ the full subcategory of $\Dsuit(\Gr(\cA))$ consisting of objects with locally coherent cohomology sheaves.
\edf

Given a dualizing complex $\omega^{\bullet}_Y \in \Dbcoh(\QCoh(Y))$, consider the functor
\begin{equation}  \label{D_Y_def}
\bD_Y \, : \, \cD_{\lc}(\Gr(\cA))^{\op} \raq \cD_{\lc}(\Gr(\cA)) \, , \qquad \cM \mapsto \RcHomcom_{\cO_Y}(\cM,\omega_Y^{\bullet})
\end{equation}
Since $\RcHomcom_{\cO_Y}(\cM,\omega_Y^{\bullet})$ can be computed in each weight component (see \eqref{RHom_weight_comp}), we have the following
\blm  \label{D_Y_inv_local}
The functor is an involution, \ie, for any $\cM \in \cD_{\lc}(\Gr(\cA))$, the canonical map $\cM \ra \bD_Y(\bD_Y(\cM))$ is an isomorphism in $\cD_{\lc}(\Gr(\cA))$. This involution interchanges $\cD_{\lc}^+(\Gr(\cA))$ and $\cD_{\lc}^-(\Gr(\cA))$.
Moreover, $\bD_Y$ is local on the base $Y$ in the sense that $\bD_Y(\cM)|_U \cong \bD_U(\cM|_U)$ for any open subscheme $U \subset Y$, where $\bD_U$ is defined using the restriction of the dualizing complex $\omega_Y^{\bullet}$ to $U$.
\elm

For each open affine subscheme $U = \Spec \, R \subset Y$, take $\scI^+(U) := \cA(U)_{>0} \cdot \cA(U) \subset \cA(U) $. It is clear that $\scI^+(U)_f = \scI^+(U_f)$ for each $f \in R$, and hence $\scI^+ \subset \cA$ forms a quasi-coherent sheaf of graded ideals of $\cA$. 
Define $\scI^-$ in a similar way.
We claim that
\begin{equation} \label{Ce_coherent}
\parbox{40em}{For each weight grading $i \in \bZ$, we have $\RGam_{\scI^+}(\cA)_i \in \Dbcoh(\QCoh(Y))$ and 
	$\Ce_{\scI^+}(\cA)_i \in \Dbcoh(\QCoh(Y))$.}
\end{equation}
Indeed, in view of Lemma \ref{RGam_affine_open}, it suffices to prove this for affine $Y = \Spec \, R$, which follows from Lemma \ref{RGam_Dbcoh}, because $\cA_0$ is assumed to be finite over $\cO_Y$ in our standing assumption \eqref{sheaf_A_setting}. 

Our discussion establishes the following

\bpp  \label{cE_RHom_comp}
There is a pair of adjunctions
\begin{equation}  \label{Ce_cE_adj_nonaffine}
\begin{tikzcd}
\cD(\Gr(\cA)) \ar[rr, shift left, "\Ce_{\scI^+}"]
& & \cD(\Gr(\cA)) \ar[rr, shift left, "(-)_0"] \ar[ll, shift left, "\cE_{\scI^+}"]
& & \cD(\QCoh(Y)) \ar[ll, shift left, "\RcHomcom_{\cO_Y}(\cA\text{,}-)"]
\end{tikzcd}
\end{equation}
where the composition $\cE_{\scI^+} \circ \RcHomcom_{\cO_Y}(\cA\text{,}-)$ is isomorphic to the functor 
$\RcHomcom_{\cO_Y}( \Ce_{\scI^+}(\cA) \text{,}-)$.

Moreover, if $\cG \in \cD^+(\QCoh(Y))$, then this composition is local on the base $Y$ in the sense that, for any affine open subscheme $U = \Spec \, R \subset Y$, if we let $A := \cA(U)$ and $G := \cG(U)$, then there is a canonical isomorphism
\begin{equation*}
\RcHomcom_{\cO_Y}( \Ce_{\scI^+}(\cA) \text{,} \cG)(U) \, \cong \, 
\RHomcom_R( \Ce_{I^+}(A) , G )
\qquad \text{ in } \cD(\GrA)
\end{equation*}
\epp

\bpf
The adjunction on the left has been shown in \eqref{RGam_LLam_Ce_cE_adj}. The one on the right is a standard adjunction. The fact that the left pointing arrows compose to the functor $\RcHomcom_{\cO_Y}( \Ce_{\scI^+}(\cA) \text{,}-)$ is an instance of \eqref{RHom_A_Y_adj}, for $\cM = \Ce_{\scI^+}(\cA)$ and $\cN = \cA$.
For the final statement, apply Lemma \ref{RHom_Y_qcoh} to $\cM = \Ce_{\scI^+}(\cA)$, which is valid in view of \eqref{Ce_coherent}.
\epf

\vspace{0.2cm}

We now discuss Serre's equivalence in the setting \eqref{sheaf_A_setting} of non-affine base. In fact, all the discussion carries through in this case almost without change. 


Let $X^+ := \Proj^+_Y(\cA)$ and $\pi^+ : X^+ \ra Y$ the projection, then there are functors
\begin{equation}  \label{GrA_QCoh_three_functors}
\begin{split}
(-)^{\sim} \,&:\, \Gr(\cA) \raq \QCoh(X^+)\, , \qquad M \, \mapsto \, \widetilde{M} \\
\!\,^0 \! \cL^{+} \,&:\, \QCoh(X^+) \raq \Gr(\cA) \, , \qquad \!\,^0 \! \cL^{+}(\cF)_i \, := \, \pi^+_*( \cF \otimes \widetilde{A(i)})
\end{split}
\end{equation}

For any $\cM\in \Gr(\cA)$, let $\!\,^0 \Ce_{\scI^+} : \Gr(\cA) \ra \Gr(\cA)$ be the functor $\!\,^0 \Ce_{\scI^+}(\cM) := \cH^0(\Ce_{\scI^+}(\cM))$. By considering its restrictions to affine open subschemes $U \subset Y$, one can alternatively describe it as 
$\cH^0(\Ce_{\scI^+}(\cM))_i \cong \pi^+_*(\widetilde{\cM(i)})$. 
As in the affine case, there are canonical maps in $\Gr(\cA)$:
\begin{equation}  \label{MCL_maps}
\begin{split}
\cM \raq \!\,^0\Ce_{\scI^+}(\cM) &\laq \!\,^0 \! \cL^{+} (\widetilde{\cM})
\end{split}
\end{equation}

In the non-affine setting, the notion of being $\tfrac{1}{d}$-Cartier still makes sense. One can either define it in a way parallel to Definition \ref{frac_Cartier}, or simply resort to the affine case:
\bdf  \label{frac_Cartier_nonaffine}
The pair $(Y,\cA)$ in \eqref{sheaf_A_setting} is said to be \emph{(positively/negatively) $\tfrac{1}{d}$-Cartier} if for some, and hence any, open affine cover $\{U_{\alpha}\}$ of $Y$, each of the $\bZ$-graded rings $\cA(U_{\alpha})$ is so, in the sense of Definition \ref{frac_Cartier}.
\edf

We have the following analogue of Propositions \ref{X_isom_Proj_2}:

\bpp  
In the situation \eqref{sheaf_of_alg_on_X}, let $\cB$ be the quasi-coherent sheaf of $\bZ$-graded algebras over $Y$ given by $\cB_i := \pi_* (\cO(i))$.
Then the sheaf of $\bN$-graded rings $\cB_{\geq 0}$ is Noetherian.
Moreover, there is a canonical isomorphism $\varphi: X \xra{\cong} \Proj^+_Y(\cB)$ over $Y$, together with a canonical isomorphism  $\bigoplus_{i \in \bZ} \varphi^*(\widetilde{\cB(i)}) \xra{\cong} \bigoplus_{i \in \bZ} \cO(i)$ of sheaves of $\bZ$-graded algebras on $X$, such that the composition
\begin{equation*}
\cB_i \raq (\pi^+_{\cB})_* (\widetilde{\cB(i)}) \, = \,  \pi_* (\varphi^*(\widetilde{\cB(i)})) \xraq{\cong} \pi_* (\cO(i)) \, =: \, \cB_i
\end{equation*}
is the identity. As a result, $\cB$ is positively $\tfrac{1}{d}$-Cartier.

Furthermore, for any quasi-coherent sheaf $\cF$ on $X$, let $\cM \in \Gr(\cB)$ be the quasi-coherent sheaf of graded $\cB$-module given by $\cM_i := \pi_*( \cF \otimes \cO(i))$, then there is a canonical isomorphism $\varphi^*( \widetilde{\cM} ) \xra{\cong} \cF$, such that the composition 
\begin{equation*}
\cM_0 \raq (\pi^+_B)_*( \widetilde{\cM}) \, = \,  \pi_* (\varphi^*(\widetilde{\cM})) \xraq{\cong}\pi_*( \cF) \, =: \, \cM_0
\end{equation*}
is the identity.
\epp

Applying this to the space $X^+ := \Proj^+_Y(\cA)$, we get the first statement of the following result; while the rest follows from Proposition \ref{MNNp}:

\bpp  \label{Proj_GrA_L_tilde}
For any $\cF \in \QCoh(X^+)$, there is a canonical isomorphism $\widetilde{\,^0\cL^{+}(\cF)} \cong \cF$.

Moreover, if we assume that $\cA$ is positively $\tfrac{1}{d}$-Cartier, then for any $\cM \in \Gr(\cA)$, we have
\begin{enumerate}
	\item for any $i \in \bZ$, the map $(\!\,^0\cL^{+} (\widetilde{\cM}))_{di} \xra{\eqref{MCL_maps}} (\!\,^0\Ce_{\scI^+}(\cM))_{di}$ is an isomorphism.
	\item both the kernel and cokernel of $\cM \xra{\eqref{MCL_maps}} \!\,^0\Ce_{\scI^+}(\cM)$ are $(\scI^+)^{\infty}$-torsion.
	\item if $\cM_{\geq c}$ is coherent over $\cA_{\geq 0}$ for some $c \in \bZ$, then so are $(\!\,^0\Ce_{\scI^+}(\cM))_{\geq c'}$ and $(\!\,^0\cL^{+} (\widetilde{\cM}))_{\geq c'}$, for any $c' \in \bZ$.
\end{enumerate}
\epp

The notions in Definition \ref{uple_torsion_def} are local on the base $Y$, and therefore extends directly to notions on $\Gr(\cA)$ in the setting \eqref{sheaf_A_setting}. This allows us to define the Serre quotients
\begin{equation*}
\QpGr(\cA) \, := \, \Gr(\cA) / \Torp(\cA)
\qquad \text{and} \qquad 
\QpmGr(\cA) \, := \, \Gr(\cA) / \Torp_{(m)}(\cA)
\end{equation*}

Proposition \ref{Proj_GrA_L_tilde} then implies the following
\bthm[Serre's equivalence]
If the pair $(Y,\cA)$ is $\tfrac{1}{d}$-Cartier, then 
the functors $(-)^{\sim}$ and $\!\,^0 \! \cL^{+}$ in \eqref{GrA_QCoh_three_functors} descend to inverse equivalences 
\begin{equation*} 
\begin{tikzcd}
\!\,^0 \! \cL^{+} \,:\, \QCoh(X^+) \ar[r, shift left]
& \QpdGr(\cA) \, : \, (-)^{\sim}  \ar[l, shift left]
\end{tikzcd}
\end{equation*}
\ethm

Now we consider the sheafified version of the results of Section \ref{GGM_sec}, which is in fact the main reason for the extended discussion above on the sheafified Hom complexes.

First, we consider the functor
\begin{equation*}
\!\,^0 \! \cR^{+} \,:\, \QCoh(X^+) \raq \Gr(\cA) \, , \qquad \!\,^0 \! \cR^{+}(\cF)_i \, := \, \pi^+_*( \cHom_{\cO_{X^+}}( \widetilde{A(-i)} , \cF ) ) 
\end{equation*}
which is again right adjoint to $(-)^{\sim}$, and comes with a canonical map
\begin{equation}  \label{CR_map}
\!\,^0\Ce_{\scI^+}(\cM) \raq \!\,^0 \! \cR^{+} (\widetilde{\cM})
\end{equation}

Corollary \ref{Cech_der_func} can be generalized to the non-affine setting to show that $\Ce_{\scI^+} : \cD(\Gr(\cA)) \ra \cD(\Gr(\cA))$ is the right derived functor of $\!\,^0 \Ce_{\scI^+}$.
Moreover, the analogue of Lemma \ref{inj_flasque} still holds.
Hence, for each $i \in \bZ$, there is a canonical isomorphism in $\cD(\QCoh(Y))$:
\begin{equation}  \label{Ce_Ip_Rpi}
\Ce_{\scI^+}(\cM)_i \, \cong \, {\bm R}\pi^+_*(\widetilde{\cM(i)})
\end{equation}
In particular, for $\cM \in \Gr(\cA)$, this gives a description of the weight components of the map $\eta_{\cM} : \cM \ra \Ce_{\scI^+}(\cM)$ in \eqref{RGam_Ce_tri}. Namely, the following diagram in $\cD(\QCoh(Y))$ commutes:
\begin{equation}
\begin{tikzcd}
\cM_i \ar[rr, "\eta_{\cM}"] \ar[d, equal] & &  \Ce_{\scI^+}(\cM)_i \ar[d, "\eqref{Ce_Ip_Rpi}", "\cong"'] \\
\cM_i \ar[r]  &  \pi^+_*(\widetilde{\cM(i)}) \ar[r]  &   {\bm R}\pi^+_*(\widetilde{\cM(i)})
\end{tikzcd}
\end{equation}
where the bottom arrows are the canonical maps.

If we take $\cR^{+} : \cD(\QCoh(X^+)) \ra \cD(\Gr(\cA))$ as the right derived functor of $\!\,^0\cR^{+}$, then there is still an adjunction at the level of derived categories:
\begin{equation} \label{m_uple_adj_derived_nonaffine}
\begin{tikzcd}
(-)^{\sim} \, : \, \cD(\Gr(\cA)) \ar[r, shift left] 
& \cD(\QCoh(X^+)) \, : \,  \cR^{+}   \ar[l, shift left] 
\end{tikzcd}
\end{equation}

From this, one can show that%
\footnote{Since $\cR^+$ is defined as a right derived functor, one obtains a map from the left hand side to the right hand side of \eqref{cR_cF_i} by comparing with the definition \eqref{GrA_QCoh_three_functors} of $\!\,^0 \! \cR^{+}$. One then shows that this map is an isomorphism by applying the adjunction isomorphism between $(-)^{\sim}$ and $\cR^+$ for the pair $K \otimes_{\cO_Y} \cA(-i) \in \cD(\Gr(\cA))$ and $\cF \in \cD(\QCoh(X^+))$, where $K \in \cD(\QCoh(Y))$.}, 
for each $i \in \bZ$, there is a canonical isomorphism in $\cD(\QCoh(Y))$:
\begin{equation}  \label{cR_cF_i}
\cR^+(\cF)_i \, \cong \, {\bm R}\pi^+_* \RcHom_{\cO_{X^+}}(\widetilde{A(-i)}, \cF)
\end{equation}

Moreover, since $\cR^+$ is a right derived functor, by comparing with \eqref{CR_map}, we have a map
\begin{equation}  \label{Cech_R_functor_nonaffine}
\Ce_{\scI^+}(\cM) \raq \cR^+(\widetilde{\cM})
\end{equation}
This map can be described in each weight component by the following commutative diagram in $\cD(\QCoh(Y))$:
\begin{equation} \label{Cech_R_functor_nonaffine_comp}
\begin{tikzcd}
\Ce_{\scI^+}(\cM)_i \ar[r] \ar[d, "\eqref{Ce_Ip_Rpi}"', "\cong"]
& \cR^+(\widetilde{\cM})_i \ar[d, "\eqref{cR_cF_i}","\cong"']\\
{\bm R}\pi^+_* ( \widetilde{M(i)})  \ar[r]
& {\bm R}\pi^+_* \RcHom_{\cO_{X^+}}(\widetilde{A(-i)}, \widetilde{M})
\end{tikzcd}
\end{equation}
where the bottom map is the canonical one.

This gives all the ingredients to give an analogue of Theorem \ref{GGM_thm}. Namely, if we replace \eqref{GGM_comp_adj_1} and \eqref{GGM_comp_adj_2}
by \eqref{Ce_cE_adj_nonaffine} and \eqref{m_uple_adj_derived_nonaffine} respectively, then we have the following

\bthm  \label{GGM_thm_nonaffine}
For each $\cG \in \cD(\QCoh(Y))$, there is a canonical isomorphism in $\cD(\Gr(\cA))$: 
\begin{equation}  \label{cR_pi_shriek_isom_sheaf}
\cR^+ ((\pi^+)^! ( \cG ))  \, \cong \,  \cE_{\scI^+}( \RcHomcom_{\cO_Y}(\cA, \cG ) ) \, \cong \, 
\RcHomcom_{\cO_Y}(\Ce_{\scI^+}(\cA), \cG )
\end{equation}
Moreover, under the isomorphisms \eqref{cR_cF_i} and \eqref{Ce_Ip_Rpi}, the weight $i$ component of this isomorphism can be described by the commutativity of the following diagram in $\cD(\QCoh(Y))$:
\begin{equation*}  
\begin{tikzcd}
\cR^+ ((\pi^+)^! ( \cG ))_i \ar[r, "\eqref{cR_pi_shriek_isom_sheaf}", "\cong"'] \ar[d, "\eqref{cR_cF_i}"', "\cong"]  &   \RcHom_{\cO_Y}(\Ce_{\scI^+}(\cA)_{-i},\cG) \ar[d, "\eqref{Ce_Ip_Rpi}", "\cong"'] \\
{\bm R}\pi^+_* \RcHom_{\cO_{X^+}}( \widetilde{\cA(-i)} , (\pi^+)^!(\cG) ) \ar[r, "\cong"] & \RcHom_{\cO_Y}( {\bm R}\pi^+_*(\widetilde{\cA(-i)}),\cG)
\end{tikzcd}
\end{equation*}
where the horizontal map in the bottom is the local adjunction isomorphism.
\ethm

\appendix

\section{Reminders on derived categories}  \label{Dcat_app}

We collect some results about derived categories. The main result of interest to us is a relation between the derived category of a quotient and the quotients of the derived categories. More precisely, Proposition \ref{Serre_SOD_1} identifies a sufficient condition where one can identify the two. In view of Proposition \ref{Groth_loc_subcat} and Corollary \ref{Groth_torsion_localizing}, this condition is satisfied in many examples.

\bdf   \label{iter_ext_def}
Given full subcategories $\cE_1 , \ldots, \cE_n$ of a triangulated category $\cD$, we denote by $ \cE_1 * \ldots * \cE_n$ the full subcategory of $\cD$ consisting of objects $X \in \cD$ with the following property:
\begin{equation}  \label{Pos_dec}
\parbox{35em}{There exists a sequence of maps $X_{n+1} \ra X_n \ra \ldots \ra X_1$ in $\cD$ such that $X_{n+1} = 0$, $X_1 = X$ and $E_i := \cone(X_{i+1} \ra X_i) \in \cE_i$ for each $1 \leq i \leq n$}.
\end{equation} 
We will often write this as
\begin{equation*}
\begin{tikzcd} [column sep = 0] 
E_n = X_n \ar[rr] & & X_{n-1} \ar[ld] \ar[rr] & &  \ldots \ar[ld] & \ldots & \ldots \ar[rr]  & & X_2 \ar[rr] \ar[ld] & & X_1 = X \ar[ld] \\
& E_{n-1} \ar[lu, "{[1]}"] && E_{n-2} \ar[lu, "{[1]}"] & & \ldots & & E_2 \ar[lu, "{[1]}"] & & E_1 \ar[lu, "{[1]}"]
\end{tikzcd}
\end{equation*}
\edf

In the special case $n=2$, the full subcategory $\cE = \cE_1 * \cE_2 $ therefore consists of objects $X \in \cD$ such that there is an exact triangle 
\begin{equation}  \label{SOD_exact_triag}
\ldots \ra E_2 \ra X \ra E_1 \ra E_2[1] \ra \ldots
\end{equation}
where $E_1 \in \cE_1$ and $E_2 \in \cE_2$.

\brm
In \cite[(1.3.9)]{BBD82}, the full subcategory $\cE_1 * \cE_2$ was denoted as $\cE_2 * \cE_1$ instead. We prefer our notation because it is consistent with the conventional notations for semi-orthogonal decomposition.
\erm

It is clear that we have $\cE_1 * \ldots * \cE_n = \cE_1 * ( \cE_2 * \ldots * \cE_n)$. In fact, by a repeated application of the octahedral axiom of a triangulated category, one has the following associativity property:
\blm[\cite{BBD82}, Lemme 1.3.10]  \label{iter_ext_assoc}
For any $1 \leq p \leq n-1$, we have \[ \cE_1* \ldots * \cE_n
\, = \,
(\cE_1 * \ldots * \cE_p ) * ( \cE_{p+1} * \ldots * \cE_n ) ) \]
\elm

\bpf
We have observed above that $\cE_1 * \ldots * \cE_n = \cE_1 * ( \cE_2 * \ldots * \cE_n)$, so that the problem reduces to showing associativity for the case $n=3$, which is \cite[Lemme 1.3.10]{BBD82}. 
\epf

\bdf
Given full subcategories $\cE_1 , \ldots, \cE_n$ of a triangulated category $\cD$, we say that the sequence $(\cE_1,\ldots,\cE_n)$ is \emph{directed} if for all $1 \leq i < j \leq n$, we have $\Hom_{\cD}(E_j,E_i) = 0$ for all $E_i \in \cE_i$ and $E_j \in \cE_j$.

In this case, we will write $\cE = \cE_1 \vec{*} \ldots \vec{*} \cE_n$ to indicate that the sequence $(\cE_1,\ldots,\cE_n)$ is directed and $\cE$ is the full subcateogry $\cE_1 * \ldots * \cE_n$.
\edf

\blm  \label{orth_rel}
Let $\cE = \cE_1 \vec{*} \cE_2$, then we have
\begin{equation*}
\begin{split}
\cE_1 \, &= \,  \{ E \in \cE \, | \, \Hom_{\cD}(E_2,E) = 0 \, \text{ for all } E_2 \in \cE_2 \} \\
\cE_2 \, &= \, \{ E \in \cE \, | \, \Hom_{\cD}(E,E_1) = 0 \, \text{ for all } E_1 \in \cE_1 \}
\end{split}
\end{equation*}
More precisely, the collections on the right hand side are precisely those that are isomorphic to objects in $\cE_1$ and $\cE_2$ respectively.
\elm

\blm
Let $\cE = \cE_1 \vec{*} \cE_2$. If both $\cE_1$ and $\cE_2$ are triangulated, then so is $\cE$.
\elm

\bpf
Closure of $\cE$ under the shift functor $[1]$ is obvious. For the closure of taking cone, suppose we are given a map $f : X \ra X'$ in $\cE$, then consider the diagram
\begin{equation*}
\begin{tikzcd}
\ldots \ar[r] & E_2 \ar[r, "i"] \ar[d, dashed] & X \ar[r, "j"] \ar[d, "f"] & E_1 \ar[r] & \ldots \\
\ldots \ar[r] & E_2' \ar[r, "i'"]& X'\ar[r, "j'"] & E_1' \ar[r] & \ldots 
\end{tikzcd}
\end{equation*}
where $E_1, E_1' \in \cE_1$ and $E_2, E_2' \in \cE_2$.
Since $j' \circ f \circ i = 0$, the dashed arrows exist, which makes the left square commute.
An application of the $3 \times 3$-lemma in a triangulated category (see, e.g., \cite[Proposition 1.1.11]{BBD82} or \cite[Lemma 2.6]{May01}) then shows that $\cone(f) \in \cE$.
\epf

For any full subcategory $\cE \subset \cD$, denote by $\langle \cE \rangle$ the smallest strictly full triangulated subcategory of $\cD$ containing $\cE$. Then we have the following

\bcor  \label{strong_orth_tri}
Suppose $\cE_1$ and $\cE_2$ are strongly orthogonal, meaning that $\Hom_{\cD}(E_2, E_1[i]) = 0$ for all $E_1 \in \cE_1$, $E_2 \in \cE_2$ and all $i \in \bZ$.  then we have 
$\langle \cE_2 \rangle \perp \langle \cE_1 \rangle$, and 
$\langle \cE_1 \vec{*} \cE_2 \rangle = \langle\cE_1 \rangle \vec{*} \langle \cE_2 \rangle$.
\ecor

Now we focus on the case when $\cE_1 \vec{*} \ldots \vec{*} \cE_n$ is the entire triangulated category $\cD$:

\bdf
A directed sequence $(\cE_1,\ldots,\cE_n)$ of full subcategories of $\cD$ is said to be a \emph{generalized semi-orthogonal decomposition} of $\cD$ if we have $\cD = \cE_1 \vec{*} \ldots \vec{*}\cE_n$.
If each $\cE_i$ is a triangulated subcategory, then the directed sequence $(\cE_1,\ldots,\cE_n)$ is said to be a \emph{semi-orthogonal decomposition} of $\cD$. We will denote a semi-orthogonal decomposition by 
$\cD = \langle \cE_1 , \ldots ,\cE_n \rangle$.
\edf

\brm
While we will focus mainly on the case of semi-orthogonal decompositions, the more general case when the subcategories $\cE_i$ are not necessarily triangulated is also important for other purposes. For example, a $t$-structure is precisely a generalized semi-orthogonal decomposition for $n=2$, such that $\cE_2[1] \subset \cE_2$, or equivalently, $\cE_1[-1] \subset \cE_1$.
\erm

\blm  \label{t_str_functorial}
Given a generalized semi-orthogonal decomposition $\cD = \langle \cE_1, \cE_2 \rangle$ such that $\cE_2[1] \subset \cE_2$, or equivalently $\cE_1[-1]\subset \cE_1$, then the assignment $X \mapsto E_1$ and $X \mapsto E_2$ in \eqref{SOD_exact_triag} are functorial, and gives a right (resp. left) adjoint to the inclusion $\cE_2 \rinto \cD$ (resp. $\cE_1 \rinto \cD$).
\elm

\bpf
Recall that, for any $Z \in \cD$, the functor $\Hom_{\cD}(-,Z) : \cD \ra \Ab$ sends exact triangles to long exact sequences. Applying this to \eqref{SOD_exact_triag} for $Z \in \cE_1$ shows that there is a canonical bijection $\Hom_{\cD}(E_1,Z) \xra{\cong} \Hom_{\cD}(X,Z)$. The shows adjunction, and hence functoriality.
\epf

\bdf
A full triangulated subcategory $\cE \subset \cD$ is said to be \emph{right admissible} (resp. \emph{left admissible}) if the inclusion functor $\cE \rinto \cD$ has a right (resp. right) adjoint. It is said to be \emph{admissible} if it is both left and right admissible. 
\edf

\blm  \label{admissible_SOD}
A full triangulated subcategory $\cE_2 \subset \cD$ is right admissible if and only if there exists a full triangulated subcategory $\cE_1 \subset \cD$ such that $\cD = \langle  \cE_1 , \cE_2  \rangle$.
Dually, a full triangulated subcategory $\cE_1 \subset \cD$ is left admissible if and only if there exists a full triangulated subcategory $\cE_2 \subset \cD$ such that $\cD = \langle  \cE_1 , \cE_2 \rangle$.
\elm

\bpf
We have already seen the implications $``\Leftarrow"$ in Lemma \ref{t_str_functorial}. Conversely, suppose that the inclusion functor $i : \cE_2 \rinto \cD$ has a right adjoint $r : \cD \ra \cE_2$, then the adjunction counit $\id \Rightarrow ri$ is an isomorphism since $i$ is fully faithful (see, e.g., \cite[Tag 07RB]{Sta}). From this, we see that $\cone[ir(X) \ra X] \in \cE_2^{\perp}$ for all $X \in \cD$, and hence there is a semi-orthogonal decomposition $\cD = \langle \cE_2^{\perp} , \cE_2 \rangle$.
\epf

Now we investigate semi-orthogonal decompositions coming from localizing subcategories of an abelian category. We start by recalling the following

\bdf
A full subcategory $\cS$ of an abelian category $\cC$ is called a \emph{Serre subcategory}%
\footnote{Often also called a \emph{dense subcategory}, for example in \cite{Pop73}} if for any short exact sequence 
$0 \ra X' \ra X \ra X'' \ra 0$ in $\cC$, $X$ is in $\cS$ if and only if both $X'$ and $X''$ are in $\cS$.
\edf

Given any Serre subcategory $\cS \subset \cC$, there is an abelian category $\cC/\cS$, together with an exact functor $\phi_{\cS}^* : \cC \ra \cC/\cS$, which is universal with respect to this property
(see, \eg, \cite[Section 4.3]{Pop73}).
For later use, we recall the construction of the category $\cC/\cS$ (see \cite[Theorem 4.3.3]{Pop73} for details).
The set of objects of $\cC/\cS$ is the same as that of $\cC$.
For any $X,Y \in \Ob(\cC/\cS) = \Ob(\cC)$, let $\mathscr{L}(X,Y)$ be the poset 
\begin{equation*}
\mathscr{L}(X,Y) \, := \, \{ \, \text{subobjects } X' \subset X \text{ and } Y' \subset Y
\text{ such that } X/X' \in \cS \text{ and } Y' \in \cS \, \}
\end{equation*}
where we define $(X',Y') \leq (X'',Y'')$ if and only if $X' \supset X''$ and $Y' \subset Y''$.
Notice that whenever $(X',Y') \leq (X'',Y'')$, there is a canonical map
$\Hom_{\cC}(X',Y/Y') \ra \Hom_{\cC}(X'',Y/Y'')$.
The Hom-sets of $\cC/\cS$ are then given by
\begin{equation}  \label{quot_Hom}
\Hom_{\cC/\cS}(X,Y) \, = \, \lim_{(X',Y') \,\in \, \mathscr{L}(X,Y)} \, \Hom_{\cC}(X',Y/Y') 
\end{equation}
The functor $\phi_{\cS}^* : \cC \ra \cC/\cS$ is defined by sending $X$ to $X$, and by taking the map $\Hom_{\cC}(X,Y) \ra \Hom_{\cC/\cS}(X,Y)$ induced by the system 
$\Hom_{\cC}(X,Y) \ra \Hom_{\cC}(X',Y/Y'')$ of canonical maps.

Another useful fact is the following
\begin{equation}  \label{quot_mon_epi}
\parbox{40em}{For any map $f : X \ra Y$ in $\cC$, its image $\phi_{\cS}^*(f) \in \Hom_{\cC/\cS}(X,Y)$ is a monomorphism (resp. an epimorphism) if and only if $\ker(f)$ (resp. $\coker(f)$) lies in $\cS$.}
\end{equation}

Now we address the question of when the functor $\phi_{\cS}^* : \cC \ra \cC/\cS$ has a left or right adjoint.

\bdf
A Serre subcategory $\cS \subset \cS$ is said to be a \emph{localizing subcategory} if the functor 
$\phi_{\cS}^* : \cC \ra \cC/\cS$ has a right adjoint $(\phi_{\cS})_* : \cC/\cS \ra \cS$.
Dually, it is said to be a \emph{colocalizing subcategory} if the functor 
$\phi_{\cS}^* : \cC \ra \cC/\cS$ has a left adjoint $(\phi_{\cS})_! : \cC/\cS \ra \cS$.
\edf

As we will see in Proposition \ref{S_closure_localizing} below, for a Serre subcategory to be (co)localizing, 
it is necessary and sufficient for objects in $\cC$ to be ``approximated'' by $\cS$-(co)closed objects in the sense of the following

\bdf  \label{S_closed_objects}
Given a Serre subcategory $\cS \subset \cC$,
an object $M \in \cC$ is said to be \emph{$\cS$-closed} if 
$\Hom_{\cC}(S,M) = \Ext^1_{\cC}(S,M) = 0$ for all $S \in \cS$.
Dually, it is said to be \emph{$\cS$-coclosed} if 
$\Hom_{\cC}(M,S) = \Ext^1_{\cC}(M,S) = 0$ for all $S \in \cS$.
\edf 

Notice that we define $\Ext^i_{\cC}(X,Y) := \Hom_{\cD(\cC)}(X,Y[i])$, so that it is well-defined, and gives rise to the standard long exact sequences, even when $\cC$ does not have enough injectives or projectives.
We now turn to the following characterization of $\cS$-closed objects, and leave the dual version to the readers.

\blm  \label{S_closed_equiv}
Given any object $M \in \cC$, the following conditions are equivalent:
\begin{enumerate}
	\item $M$ is $\cS$-closed.
	\item Given $f : X \ra Y$ such that both $\ker(f)$ and $\coker(f)$ are in $\cS$, the induced map 
	$\Hom_{\cC}(Y,M) \ra \Hom_{\cC}(X,M)$ is a bijection.
	\item Given an injection $f : X \ra Y$ such that $\coker(f)$ is in $\cS$, the induced map 
	$\Hom_{\cC}(Y,M) \ra \Hom_{\cC}(X,M)$ is a bijection.
	\item For any $X \in \cC$, the induced map $\Hom_{\cC}(X,M) \ra \Hom_{\cC/\cS}(X,M)$
	is a bijection.
\end{enumerate}
\elm

\bpf
For $(1) \Rightarrow (2)$, apply the long exact sequence associated to $\Ext^{\bullet}_{\cC}(-,M)$.
The implication $(2) \Rightarrow (3)$ is obvious.
For $(3) \Rightarrow (1)$, given any $S \in \cS$, applying condition $(3)$ to the injection $0 \ra S$ shows that $\Hom_{\cC}(S,M) = 0$. Also, if $\Ext^1_{\cC}(S,M) \neq 0$, then there exists an extension 
$0 \ra M \ra X \ra S \ra 0$ such that the connecting homomorphism $\Hom_{\cC}(M,M) \ra \Ext^1_{\cC}(S,M)$ sends the identity to a nonzero element. This contradicts the requirement in condition $(3)$ that $\Hom_{\cC}(X,M) \ra \Hom_{\cC}(M,M)$ be a bijection.

For $(3) \Rightarrow (4)$, given any $S \in \cS$, observe as above that $\Hom_{\cC}(S,M) = 0$. 
Thus, $M$ has no subobject in $\cS$, and the description \eqref{quot_Hom} of the Hom-sets of $\cC/\cS$ for $Y = M$ therefore becomes an inverse limit over $(X',0) \in \mathscr{L}(X,M)$. Condition $(3)$ then implies that each of the canonical maps $\Hom_{\cC}(X,M) \ra \Hom_{\cC}(X',M)$ is a bijection, hence showing $(4)$.
For $(4) \Rightarrow (2)$, simply notice that a map $f : X \ra Y$ as in $(2)$ descends to an isomorphism in $\cC/\cS$, by \eqref{quot_mon_epi}.
\epf

Now if $\cS \subset \cC$ is a localizing subcategory, then for any  $M = (\phi_{\cS})_* (\cZ)$ and any $X \in \cC$, we have by adjunction
$\Hom_{\cC}(X,M) = \Hom_{\cC/\cS}(\phi_{\cS}^*(X), \cZ)$,
so that the value of the functor $\Hom_{\cC}(-,M)$ on $X$ depends only on $\phi_{\cS}^*(X)$.
Thus, if we are given $f : X \ra Y$
as in condition $(2)$ of Lemma \ref{S_closed_equiv},
then by \eqref{quot_mon_epi}, $f$ descends to an isomorphism in $\cC/\cS$, so that 
$\Hom_{\cC}(Y,M) \ra \Hom_{\cC}(X,M)$ is a bijection, and $M$ is therefore $\cS$-closed.
%
It then follows formally (see \cite[Proposition 4.4.3]{Pop73}) that the adjunction counit
$\phi_{\cS}^* (\phi_{\cS})_*(\cZ) \ra \cZ$ is an isomorphism in $\cC/\cS$.
As a result, for any $X \in \cC$, the adjunction unit $\epsilon_X : X \ra (\phi_{\cS})_* \phi_{\cS}^* (X)$ 
has $\ker(\epsilon_X) \in \cS$ and  $\coker(\epsilon_X) \in \cS$.
Thus, the adjunction unit consistutes an example of an $\cS$-closure, in the sense of the following

\bdf  \label{S_closure_def}
Given a Serre subcategory $\cS \subset \cC$, an \emph{$\cS$-closure} of an object $X \in \cC$ is a map $\epsilon_X : X \ra \bar{X}$ from $X$ to an $\cS$-closed object $\bar{X}$, such that $\ker(\epsilon_X)$ and  $\coker(\epsilon_X)$ are both in $\cS$.
\edf

Notice that, by condition $(2)$ of Lemma \ref{S_closed_equiv}, if an $\cS$-closure exists, then it is unique up to canonical isomorphism. 
This notion is useful because of the following

\bpp  \label{S_closure_localizing}
Given a Serre subcategory $\cS \subset \cC$, the following statements are equivalent:
\begin{enumerate}
	\item The Serre subcategory $\cS \subset \cC$ is a localizing subcategory.
	\item Every object in $\cC$ has an $\cS$-closure.
	\item Every object $X \in \cC$ has a largest subobject $X_{\cS}$ in $\cS$, and $X/X_{\cS}$ embeds into an $\cS$-closed object.
\end{enumerate}
\epp

\bpf
We have already seen the implication $(1) \Rightarrow (2)$ in the paragraph preceding Definition \ref{S_closure_def}.

For $(2) \Rightarrow (1)$, notice that if $Y \ra \bar{Y}$ is an $\cS$-closure, then by Lemma \ref{S_closed_equiv}(4), we have $\Hom_{\cC}(X,\bar{Y}) \cong \Hom_{\cC/\cS}(\phi^*(X) , \phi^*(Y))$ for all $X \in \cC$. Thus, the object $\bar{Y} \in \cC$ is the would-be-image of $\phi^*(Y) \in \cC/\cS$ under a would-be-right-adjoint to $\phi^*$. By Yoneda lemma, this in turn guarantees that the functor $\phi_* : \cC/\cS \ra \cC$ given by $\phi^*(Y) \mapsto \bar{Y}$ is well-defined and right adjoint to $\phi^*$.

For $(2) \Rightarrow (3)$, Suppose $\epsilon_X : X \ra \bar{X}$ is an $\cS$-closure. Then the kernel $X_{\cS} := \ker(\epsilon_X)$
is in $\cS$. Moreover, any subobject $X' \subset X$ in $\cS$ must map to zero in $\bar{X}$ since $\Hom_{\cC}(X',\bar{X}) = 0$ by $\cS$-closedness of $\bar{X}$. Thus, $X_{\cS}$ is the largest subobject of $X$ in $\cS$, and we have an emebedding $X / X_{\cS} \rinto \bar{X}$.

For $(3) \Rightarrow (2)$, suppose we are given a subobject $X_{\cS}$ of $X$ in $\cS$, and an embedding $i : X / X_{\cS} \rinto M$ into
an $\cS$-closed object $M$. Let $K := \coker(i)$, with the canonical epimorphism $j : M \ronto K$,
and let $\bar{X} = j^{-1}(K_{\cS})$. Thus we have an exact sequence 
$0 \ra X_{\cS} \ra X \ra \bar{X} \ra K_{\cS} \ra 0$, with $X_{\cS}$ and $K_{\cS}$ both in $\cS$.
Therefore, it suffices to show that $\bar{X}$ is $\cS$-closed.
To this end, observe first that we have $j : M/\bar{X} \xra{\cong} K/K_{\cS}$. Since $K/K_{\cS}$ has no subobjects in $\cS$, this implies that $\Hom_{\cC}(S , M/\bar{X}) = 0$ for all $S \in \cS$.
The $\cS$-closedness of $\bar{X}$ then follows by applying, for each $S \in \cS$, the $\Ext^{\bullet}(S,-)$ long exact sequence 
to the short exact sequence $0 \ra \bar{X} \ra M \ra M/\bar{X} \ra 0$,
using the $\cS$-closedness of $M$.
\epf

We record the following simple Lemma for later use:

\blm  \label{lift_SES}
Let $\cS \subset \cC$ be a localizing subcategory, then short exact sequences in $\cC/\cS$
can be functorially lifted to short exact sequences in $\cC$.
\elm

\bpf
Given $0 \ra \cX \xra{f} \cY \xra{g} \cZ \ra 0$, take the short exact sequence $0 \ra \phi_* \cX \xra{\phi_*(f)} \phi_* \cY \ra \coker(\phi_*(f)) \ra 0$. The functor $\phi^*$ then sends it to the one we start with because it is exact.
\epf

We now identify a situation when a given Serre subcategory is easily recognized to be localizing. 
See Proposition \ref{left_reg_tor_prop} and Corollary \ref{Groth_torsion_localizing} below.

\bdf  \label{tor_def}
A \emph{torsion pair} for an abelian category $\cC$ consists of a pair of full subcategories $\cF$ and $\cT$ satisfying the following conditions:
\begin{enumerate}
	\item $\cF \cap \cT = 0$ ;
	\item If $X$ is an object of $\cT$, then any quotient object of $X$ is also in $\cT$;
	\item If $X$ is an object of $\cF$, then any subobject of $X$ is also in $\cF$;
	\item For each $X \in \cC$, there is an exact sequence
	$
	0 \ra X^{\cT} \ra X \ra X_{\cF} \ra 0
	$
	with $X^{\cT} \in \cT$ and $X_{\cF} \in \cF$.
\end{enumerate}
We write $\cC = \langle \cF,\cT \rangle$ for a torsion pair.
\edf

\blm  \label{orth_rel_tor}
Let $\cC = \langle \cF, \cT \rangle$ be a torsion pair, then we have
\begin{equation*}
\begin{split}
\cF \, &= \, \{ X \in \cC \, | \, \Hom_{\cC}(T,X) = 0 \, \text{ for all } T \in \cT \} \\
\cT \, &= \, \{ X \in \cC \, | \, \Hom_{\cC}(X,F) = 0 \, \text{ for all } F \in \cF \}
\end{split}
\end{equation*}
More precisely, the collections on the right hand side are precisely those that are isomorphic to objects in $\cF$ and $\cT$ respectively.
\elm

\bdf  \label{left_reg_tor}
A torsion pair $\cC = \langle \cF,\cT \rangle$ is said to be \emph{injectively cogenerated} if the following two conditions hold:
\begin{enumerate}
	\item If $X$ is an object of $\cT$, then any subobject of $X$ is also in $\cT$;
	\item For any $F \in \cF$, there is a monomorphism $F \rinto I$ where $I$ is an injective object of $\cC$ that lies in $\cF$.
\end{enumerate}
\edf

Given a torsion pair $\cC = \langle \cF,\cT \rangle$ satisfying condition (1) of this Definition, then $\cT$ is a Serre subcategory. If $I$ is an injective object of $\cC$ that lies in $\cF$, then by Lemma \ref{orth_rel_tor}, the object $I$ is $\cT$-closed in the sense of Definition \ref{S_closed_objects}. Thus, by Proposition \ref{S_closure_localizing}(3), we have (the first statement of) the following

\bpp  \label{left_reg_tor_prop}
If $\cC = \langle \cF,\cT \rangle$ is an injectively cogenerated torsion pair, then 
\begin{enumerate}
	\item  $\cT \subset \cC$ is a localizing subcategory;
	\item If $I$ is an injective object of $\cC$ that lies in $\cF$, then $\phi^*(I)$ is injective in $\cC/\cT$;
	\item The category $\cC/\cT$ has enough injectives.
\end{enumerate}
\epp

\bpf
In the paragraph preceding this Proposition, we have already shown (1), as well as the fact that the object $I$ in (2) is $\cT$-closed. Thus, by Lemma \ref{S_closed_equiv}(4), we have $\Hom_{\cC/\cT}(X,I) \cong \Hom_{\cC}(\phi^*(X), \phi^*(I))$ for all $X \in \cC$. By Lemma \ref{lift_SES}, epimorphisms in $\cC/\cS$ can be functorially lifted to epimorphisms in $\cC$. It follows that $\phi^*(I)$ is injective in $\cC/\cT$.

To show (3), given $\phi^*(X) \in \cC/\cT$, choose a monomorphism $X_{\cF} \rinto I$ as in Definition \ref{left_reg_tor}. Thus, we have an exact sequence $0 \ra X^{\cT} \ra X \ra I$. Applying the exact functor $\phi^*$ gives $0 \ra \phi^*(X) \ra \phi^*(I)$.
\epf

\blm
If $\cC = \langle \cF,\cT \rangle$ is a torsion pair such that $\cT$ is closed under subobjects (\ie, if it satisfies condition (1) of Definition \ref{left_reg_tor}) then any essential extension of an object $F \in \cF$ is still in $\cF$.
\elm

\bpf
Given an essential extension $F \rinto X$, if there exists some nonzero map $f : T \ra X$ for some $T \in \cT$, then its image ${\rm Im}(f)$ is a nonzero subobject of $X$ that lies in $\cT$. By assumption, ${\rm Im}(f)$ has a nonzero intersection $Z$ with $F$. Since $Z \subset F$, we have $Z \in \cF$. Also, by the assumption that $\cT$ is closed under subobjects, we have $Z \in \cT$. Thus, $Z = 0$, a contradiction.
\epf

\bcor  \label{inj_env_reg}
If $\cC$ admits injective envelope, then every torsion pair in $\cC$ that satisfies condition (1) of Definition \ref{left_reg_tor} also satisfies condition (2).
\ecor

This is useful because of the following well-known result (see, e.g., \cite[Theorem 3.10.10]{Pop73}):
\bthm  \label{Groth_inj_env}
Every Grothendieck category admits injective envelopes.
\ethm

Combining Corollary \ref{inj_env_reg}, Theorem \ref{Groth_inj_env} and Proposition \ref{left_reg_tor_prop}, we have the following 

\bcor  \label{Groth_torsion_localizing}
Let $\cC$ be a Grothendieck category. If
$\cC = \langle \cF,\cT \rangle$  is a torsion pair such that $\cT$ is closed under subobjects (\ie, if condition (1) of Definition \ref{left_reg_tor} is satisfied), then the Serre subcategory $\cT$ is localizing.
\ecor

Now, we identify a sufficient criterion where one can identify the derived category of a quotient with the quotients of the derived categories. See Proposition \ref{Serre_SOD_1} below.

Let $\cC$ be an abelian category. Denote by $\cD(\cC)$ its derived category.
For any Serre subcategory
$\cS \subset \cC$, denote by $\cD_{\cS}(\cC) \subset \cD(\cC)$ the full subcategory consisting of complexes whose cohomology lie in $\cS$.
Then $\cD_{\cS}(\cC)$ is a split-closed triangulated subcategory%
\footnote{This holds more generally for weak Serre subcategory in the sense of \cite[Tag 02MO]{Sta} (see \cite[Tag 06UQ]{Sta})}  of $\cD(\cC)$.
In fact, since the canonical functor $\phi^* : \cC \ra \cC/\cS$ is exact, it descends to an exact functor 
$
\phi^* : \cD(\cC) \ra \cD(\cC/\cS)
$, 
whose kernel is precisely $\cD_{\cS}(\cC)$.
Similar statements hold when $\cD$ is replaced by $\cD^+$, $\cD^-$ or $\cD^b$.



\bdf  \label{D_localizing_def}
A Serre subcategory $\cS \subset \cC$ is said to be \emph{$\Dsuit$-localizing} (where $\spadesuit \in \{\, \, \,   ,+,-,b\}$) if the functor $\phi^* : \Dsuit(\cC) \ra \Dsuit(\cC/\cS)$ has a right adjoint ${\bm R}\phi_* : \Dsuit(\cC/\cS) \ra \Dsuit(\cC)$ such that the adjunction counit $ \epsilon : \phi^* \circ {\bm R}\phi_* \Rightarrow \id$ is an isomorphism.


Dually, $\cS$ is said to be \emph{$\Dsuit$-colocalizing} if the functor $\phi^* : \Dsuit(\cC) \ra \Dsuit(\cC/\cS)$ has a left adjoint ${\bm L}\phi_! : \Dsuit(\cC/\cS) \ra \Dsuit(\cC)$
such that the adjunction unit $\eta :  \id \Rightarrow \phi^* \circ {\bm L}\phi_! $ is an isomorphism.
\edf


Although we do not assume in the above Definition that ${\bm R}\phi_*$ and ${\bm L}\phi_!$ are derived functors of some functors at the abelian level, this will often be the case in applications. The following result gives the main class of example of $\cD$-localizing subcategory:

\bpp \label{Groth_loc_subcat}
Let $\cC$ be a Grothendieck abelian category and $\cS \subset \cC$ any localizing subcategory. Then we have
\begin{enumerate}
	\item The Serre quotient $\cC/\cS$ is a Grothendieck category.
	\item For $\spadesuit \in \{\, \, \,  ,+\}$, the functor $\phi_* : \cC/\cS \ra \cC$ has a total right derived functor ${\bm R}\phi_* : \Dsuit(\cC/\cS) \ra \Dsuit(\cC)$, which is right adjoint to $\phi^* : \Dsuit(\cC) \ra \Dsuit(\cC/\cS)$, and makes $\cS \subset \cC$ a $\Dsuit$-localizing subcategory.
\end{enumerate}
\epp

\bpf
The quotient functor $\phi^* : \cC \ra \cC/\cS$ preserves arbitrary colimits since it has a right adjoint. Thus $\cC/\cS$ admits small colimits. Since short exact sequences in $\cC/\cS$ can be functorially lifted to short exact sequences in $\cC$ (see Lemma \ref{lift_SES}), directed colimit is exact in $\cC/\cS$. Moreover, the functor $\phi^*$ sends any generating set of $\cC$ to a generating set of $\cC/\cS$. Thus, $\cC/\cS$ is a Grothendieck category, proving (1).

Since any Grothendieck category has enough injectives (see, e.g., \cite[Tag 079H]{Sta}) and K-injectives (see, e.g., \cite[Tag 079P]{Sta}), the functor $\phi_* : \cC/\cS \ra \cC$ can be derived to ${\bm R}\phi_* : \Dsuit(\cC/\cS) \ra \Dsuit(\cC)$ for $\spadesuit \in \{\, \, \,  ,+\}$.
Moreover, as a right adjoint to an exact functor, the functor $\phi_* : \cC/\cS \ra \cC$ preserves injectives and K-injectives. Thus, the derived functor ${\bm R}\phi_* : \Dsuit(\cC/\cS) \ra \Dsuit(\cC)$ is right adjoint to $\phi^* :  \Dsuit(\cC) \ra  \Dsuit(\cC/\cS)$. The fact that the adjunction counit $ \epsilon : \phi^* \circ {\bm R}\phi_* \Rightarrow \id$ is an isomorphism is also clear by applying it on any (K-)injective representative.
\epf

The usefulness of Definition \ref{D_localizing_def} lies in the following obvious

\bpp   \label{Serre_SOD_1}
Suppose that $\cS \subset \cC$ is a $\Dsuit$-localizing subcategory, 
then ${\bm R}\phi_*$ is fully faithful, and gives rise to a semi-orthogonal decomposition 
\begin{equation*}
\Dsuit(\cC) \, = \, \langle \,  {\bm R}\phi_* ( \Dsuit(\cC/\cS)) \,   , \, \Dsuit_{\cS}(\cC) \, \rangle
\end{equation*}

As a result, there is an equivalence of triangulated categories
\begin{equation*}
\begin{tikzcd}
\phi^* \, : \, \Dsuit(\cC) / \Dsuit_{\cS}(\cC) \ar[r, shift left] & \Dsuit(\cC/\cS) \, : \, {\bm R} \phi_* \ar[l, shift left]
\end{tikzcd}
\end{equation*}
\epp

\brm
An analogue of Proposition \ref{Serre_SOD_1} is claimed in \cite[Tag 06XM]{Sta}. However, the proof seems to be incomplete.
\erm

We now recall the notion of (sequential) homotopy limits and (sequential) homotopy colimits of complexes.
Let $\cC$ be an abelian category with exact small coproducts (also known as an Ab 4-category),
then for any directed system $X_{\bullet} = [X_0 \xra{f_0} X_1 \xra{f_1} \ldots]$ of objects in $\Ch(\cC)$, define its \emph{homotopy colimt} to be the cochain complex defined by the cone
\begin{equation}  \label{hocolim_def}
\hocolim_{n \in \bN} X_{n} \, := \, \cone \, [ \,  \coprod_{n \in \bN} X_n \xraq{\alpha} \coprod_{n \in \bN} X_n   \, ] 
\end{equation}
where the map $\alpha$ is the coproduct of the maps $X_n \xra{\id_{X_n} -f_n} X_n \amalg X_{n+1}$.

Dually, if $\cC$ be an abelian category with exact small products (also known as an Ab 4*-category),
then for any inverse system $X^{\bullet} = [X^0 \xla{f^0} X^1 \xla{f^1} \ldots]$ of objects in $\Ch(\cC)$, define its \emph{homotopy limt} to be the cochain complex defined by the cocone
\begin{equation}   \label{holim_def}
\holim_{n \in \bN} X^{n} \, := \, \cocone \, [ \,  \prod_{n \in \bN} X^n \xraq{\beta} \prod_{n \in \bN} X^n   \, ] 
\end{equation}
where the map $\beta$ is the product of the maps $X^n \oplus X^{n+1} \xra{(\id_{X^n}, -f^n)} X^n$.

In general, the ordinary (termwise) directed colimit of a directed system $X_{\bullet}$ in $\Ch(\cC)$ is the cokernel of the map $\alpha$ in \eqref{hocolim_def}; while the ordinary (termwise) inverse limit of an inverse system $X^{\bullet}$ in $\Ch(\cC)$ is the kernel of the map $\beta$ in \eqref{holim_def}. Thus, we have canonical maps
\begin{equation}  \label{hocolim_to_colim}
\begin{split}
\hocolim_{n \in \bN} X_{n} \raq \varinjlim X_n  \qquad &\text{ for direct system } X_{\bullet} \\
\varprojlim X^n \raq \holim_{n \in \bN} X^{n} \qquad &\text{ for inverse system } X^{\bullet}
\end{split}
\end{equation}

The notion of homotopy colimits often coincides with the ordinary colimits, in view of the following Lemma (see, e.g., \cite[Tag 0949]{Sta}):

\blm  \label{hocolim_colim_qism}
If $\cC$ is an Ab 5 category (\ie, directed colimits exist and are exact), then for any directed system $X_{\bullet}$, the canonical map $\hocolim_{n \in \bN} X_{n} \raq \varinjlim X_n$ is a quasi-isomorphism.
\elm

Finally, we record the following well-known result:

\blm  \label{Dbcoh_pc_lem}
Let $\cC$ be an abelian category, and let $\cP \subset \Ob(\cC)$ be a collection of projective objects closed under finite direct sum. Denote by $Q(\cP) \subset \Ob(\cC)$ the collection of objects $M$ such that there exists an epimorphism $P \ronto M$ from some $P \in \cP$. 
Suppose $Q(\cP)$ is closed under taking subobjects, then for any bounded above complex $M^{\bullet}$ in $\cC$ whose cohomology objects lie in $Q(\cP)$, there exists a bounded above complex $P^{\bullet}$ of objects in $\cP$, together with a quasi-isomorphism $\varphi: P^{\bullet} \rsa M^{\bullet}$.
\elm

\bpf
Assume that $M^i = 0$ for $i > b$. Choose an epimorphism $P^b \ronto M^b / d(M^{b-1})$, and lift it to a map $P^b \ra M^b$. Suppose there is a complex $P^{\bullet}$ of objects of $\cP$ concentrated in degrees $[a,b]$, together with a map $\varphi: P^{\bullet}\ra M^{\bullet}$, such that $H^i(P^{\bullet}) \ra H^{i}(M^{\bullet})$ is an isomorphism for $i > a$ and is surjective for $i=a$. 
The last condition guarantees that the map 
\begin{equation}  \label{phi_d_P_M_complex}
(\varphi , -d) \, : \, Z^a(P^{\bullet}) \oplus (M^{a-1}/d(M^{a-2})) \raq Z^a(M^{\bullet})
\end{equation}
is surjective. Notice that $Q(\cP)$ is a Serre subcategory, so that $\cD_{Q(\cP)}(\cC)$ is a triangulated subcategory. Since the two term complex $[ M^{a-1}/d(M^{a-2}) \xra{-d} Z^a(M^{\bullet})]$
has cohomology objects $H^{a-1}(M^{\bullet})$ and  $H^{a}(M^{\bullet})$, and since the map \eqref{phi_d_P_M_complex}, thought of as a two term complex, is an extension of it by $Z^a(P^{\bullet}) \in Q(\cP)$, we see that the kernel $K$ of \eqref{phi_d_P_M_complex} is in $Q(\cP)$ as well. Choose an epimorphism $P' \ronto K$ with $P' \in \cP$. Then we have a commutative diagram
\begin{equation*}
\begin{tikzcd}
P' \ar[r, "d'"] \ar[d, "\psi'"'] & Z^a(P^{\bullet}) \ar[d, "\varphi"] \\
M^{a-1}/d(M^{a-2}) \ar[r, "d"] & Z^a(M^{\bullet}) 
\end{tikzcd}
\end{equation*}
such that the induced map $\varphi : Z^a(P^{\bullet}) / d'(P') \ra Z^a(M^{\bullet}) / d(M^{a-1}) = H^a(M^{\bullet})$ is an isomorphism. 
Moreover, choose an epimorphism $\psi'' : P'' \ronto H^{a-1}(M^{\bullet}) \subset M^{a-1}/d(M^{a-2})$
from some $P'' \in \cP$, and consider the commutative diagram 
\begin{equation*}
\begin{tikzcd}
P' \oplus P'' \ar[r, "(d'\text{,}0)"] \ar[d, "(\psi' \text{,} \psi'')"'] & Z^a(P^{\bullet}) \ar[d, "\varphi"] \\
M^{a-1}/d(M^{a-2}) \ar[r, "d"] & Z^a(M^{\bullet}) 
\end{tikzcd}
\end{equation*}
Let $P^{a-1} := P' \oplus P''$, and lift the map $(\psi' , \psi'')$ to a map $\varphi : P^{a-1} \ra M^{a-1}$.
This completes the desired inductive step.
\epf


\begin{thebibliography}{9}
\bibitem{AZ94}
M. Artin, and J. Zhang,
\textit{Noncommutative projective schemes}, 
Adv. Math. \textbf{109} (1994), 228--287. 

\bibitem{AKO08}
D. Auroux, L. Katzarkov, and D. Orlov, 
\textit{Mirror symmetry for weighted projective planes and their noncommutative deformations}, 
Ann. of Math. (2) \textbf{167} (2008), 867--943. 

\bibitem{ATJLJ97}
L. Alonso Tarrío, A. Jerem\'{i}as L\'{o}pez, and J. Lipman,
\textit{Local homology and cohomology on schemes}, 
Ann. Sci. \'{E}cole Norm. Sup. (4) \textbf{30} (1997), 1--39. 




\bibitem{BH93}
W. Bruns, and J. Herzog, 
\textit{Cohen-Macaulay rings}, 
Cambridge Studies in Advanced Mathematics, \textbf{39}. Cambridge University Press, Cambridge, 1993.

\bibitem{BBD82}
A. Beĭlinson, J. Bernstein, and P. Deligne,
\textit{Faisceaux pervers}, 
Analysis and topology on singular spaces, I (Luminy, 1981), 5--171, 
Astérisque, 100, Soc. Math. France, Paris, 1982. 




















\bibitem{HL15}
D. Halpern-Leistner, 
\textit{The derived category of a GIT quotient},
J. Amer. Math. Soc. \textbf{28} (2015), 871--912. 









\bibitem{KM98}
J. Koll\'{a}r, and S. Mori, \textit{Birational geometry of algebraic varieties}, With the collaboration of C. H. Clemens and A. Corti. Translated from the 1998 Japanese original. Cambridge Tracts in Mathematics, 134. Cambridge University Press, Cambridge, 1998.










\bibitem{May01}
J. P. May, 
\textit{The additivity of traces in triangulated categories},
Adv. Math. \textbf{163} (2001), no. 1, 34--73.







\bibitem{Pop73}
N. Popescu, 
\textit{Abelian categories with applications to rings and modules}, 
London Mathematical Society Monographs, Academic Press, London-New York, 1973.








\bibitem{Sta}
The Stacks Project Authors,
\textit{Stacks Project},
\url{https://stacks.math.columbia.edu},
2019



\bibitem{Tha96}
M. Thaddeus,
\textit{Geometric invariant theory and flips}, 
J. Amer. Math. Soc. \textbf{9} (1996), no. 3, 691--723. 

\bibitem{TT90}
R. Thomason, and T. Trobaugh, \textit{Higher algebraic K-theory of schemes and of derived categories}, The Grothendieck Festschrift, Vol. III, 247–435, Progr. Math., \textbf{88}, Birkh\"{a}user Boston, Boston, MA, 1990.













\bibitem{Yeu20a}
W.K. Yeung, 
\textit{Homological flips and homological flops}, {\tt arXiv:1907.06190}

\bibitem{Yeu20b}
W.K. Yeung, 
\textit{Weight truncation for wall-crossings in birational cobordisms}, preprint
















\end{thebibliography}
\end{document}